\documentclass[a4paper, 12pt]{amsart}
\usepackage[T2A]{fontenc}
\usepackage[utf8]{inputenc}
\usepackage[english]{babel}

\usepackage{amsmath,amssymb,amsthm,url}
\usepackage{dsfont,textcomp,graphicx,overpic}
\usepackage[usenames,dvipsnames]{color}

\graphicspath{{figs/}}

\makeatletter

\@addtoreset{figure}{section}
\makeatother

\usepackage{hyperref}
\hypersetup{
   colorlinks   = true, %Colours links instead of ugly boxes
   urlcolor     = blue, %Colour for external hyperlinks
   linkcolor    = blue, %Colour of internal links
   citecolor   = red %Colour of citations
}

\usepackage{tikz}
\usetikzlibrary{positioning, calc, intersections, through, arrows}
\usetikzlibrary{graphs}
\usetikzlibrary{graphs.standard}
\usetikzlibrary{patterns}

\def\sgn{\mathop{\fam0 sgn}}

\def\wu{\mathop{\fam0 wu}}

\def\R{{\mathbb R}} \def\Z{{\mathbb Z}}  

\theoremstyle{plain}
\newtheorem{theorem}{Theorem}[section]

    \newtheorem{example}[theorem]{Example}

\newtheoremstyle{mydefinition}% name
  {\medskipamount}%      Space above
  {\medskipamount}%      Space below
  {\normalfont}%         Body font
  {\parindent}%         Indent amount (empty = no indent, \parindent = para indent)
  {\bfseries}% Thm head font
  {.}%        Punctuation after thm head
  { }%     Space after thm head: " " = normal interword space;
  {}%         Thm head spec (can be left empty, meaning `normal')

\theoremstyle{mydefinition}
\newtheorem{pr}[theorem]{}

\begin{document}

\title{Invariants of almost embeddings of graphs in the plane: results and problems} 
 
\author{E. Alkin, E. Bordacheva, A. Miroshnikov, O. Nikitenko, A. Skopenkov}
\thanks{\emph{E. Alkin, E. Bordacheva, A. Miroshnikov, A. Skopenkov:} Moscow Institute of Physics and Technology. 
\newline
\emph{O. Nikitenko:} Altay Technical University. 
\newline
\emph{A. Skopenkov:} Independent University of Moscow, \url{https://users.mccme.ru/skopenko/}. 
\newline
We are grateful to M. Didin, P. Kozhevnikov, and especially to T. Garaev for useful discussions. 
\newline
MSC-class: 55-02, 55M25, 05C10, 57K20. 
}
 
\date{}

\maketitle

\begin{abstract} 
A graph drawing in the plane is called an \textit{almost embedding} if images of any two non-adjacent simplices (i.e. vertices or edges) are disjoint.
We introduce integer invariants of almost embeddings: winding number, 
%(\S\ref{s:wind}), 
cyclic and triodic Wu numbers. 
%(\S\ref{s:gawh}). 
%(In this text an invariant is just a number assigned to an almost embedding, we do not introduce the notion of almost isotopy under which this number is invariant.)
%Which values these invariants can assume? 
We construct almost embeddings realizing some values of these invariants. 
We prove some relations between the invariants.  
We study values realizable as invariants of some almost embedding, but not of any embedding. 
%More precisely, main results are assertions \ref{borul-pl}, \ref{t:radonae}, \ref{e:radonae}, \ref{t:k5-e}, \ref{p:triod}.b, \ref{p:off}.b, \ref{e:k6}.a.   

This paper is expository and is accessible to mathematicians not specialized in the area (and to students). 
However elementary, this paper is motivated by frontline of research. 
\end{abstract}

\tableofcontents

%\newpage
%\section*{Preface}
  
%You can postpone reading this section, and start solving problems from \S\ref{s:wind}. 
%\textbf{Informal description of the main subject and main results}

\begin{figure}[ht]\centering
        \includegraphics[width=4.6cm]{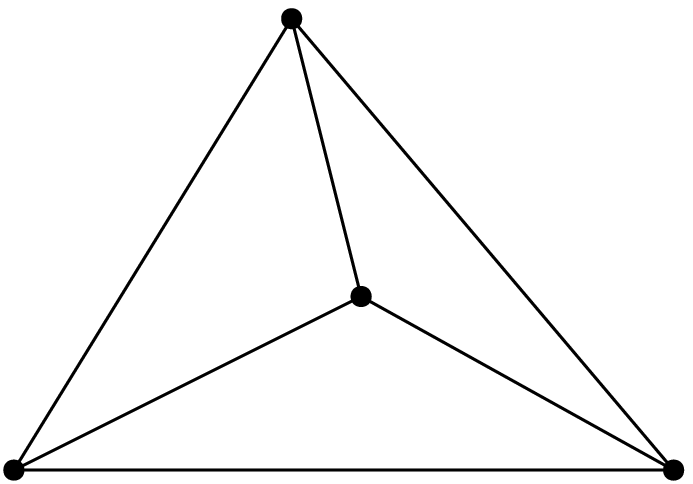} \qquad 
        \includegraphics[width=4.6cm]{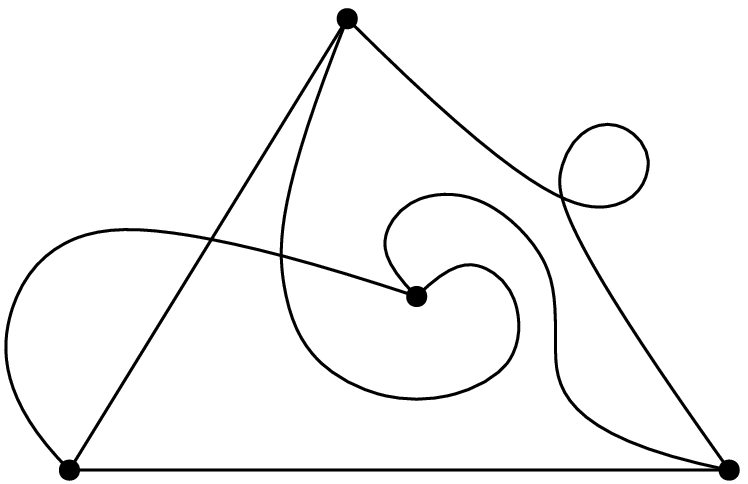} \qquad
        \includegraphics[width=3cm]{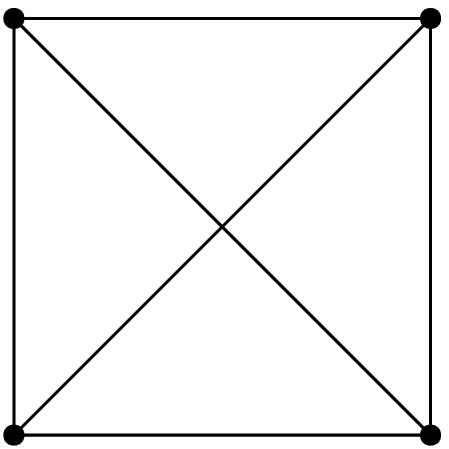}
    \caption{An embedding, an almost embedding, and a map (drawing) which is not an almost embedding}\label{f:k4} 
    \end{figure}

A classical subject is study of planar graph drawings without self-intersections (i.e. embeddings or plane graphs).   
It is also interesting to study graph drawings having `moderate' self-intersections, e.g.  almost embeddings 
%is (roughly speaking) a graph drawing for which the images of any two non-adjacent edges are disjoint, 
(see definition near Figures \ref{f:k4} and \ref{k5}).   

For relation to modern research see remark at the end of \S\ref{s:alem},  
%\ref{r:mot}, 
papers \cite{Sk18, KS20, IKN+, Ga23, Bo} and the references therein.
We do not know proofs of Conjectures \ref{c:radonae}.b, \ref{p:k33}.b, \ref{p:wu-conj}.b, \ref{p:wi-wu}.b, \ref{e:k6}.b, or solution of Problem \ref{p:cube}. 

%\medskip
%\textbf{Learning by doing problems}
%Here complicated material is split into simple steps and is illustrated by examples (this is one of the reasons for the length of this text).   
 
\emph{In this text we expose a theory as a sequence of problems,} see e.g. \cite{HC19}, \cite[Introduction, Learning by doing problems]{Sk21m} and the references therein. 
Most problems are useful theoretical facts. 
So this text could in principle be read even without solving problems or looking to \S\ref{s:answ}. 
Problems are numbered, the words `problem' are omitted. 
If a mathematical statement is formulated as a problem, then the objective is to prove this statement.
Open-ended questions are called {\bf riddles}; here one must come up with a clear wording, and a proof. 
\emph{If a problem is named `theorem' (`lemma', `corollary', etc.), then this statement is considered to be more important.}
Usually we \emph{formulate} beautiful or important statement \emph{before} giving a sequence of results (lemmas, assertions, etc.) which constitute its \emph{proof}.
We give hints on that after the statements but we do not want to deprive you of the pleasure of finding the right moment when you finally are ready to prove the statement.
%In general, if you are stuck on a certain problem, try looking at the next ones; they may turn out to be helpful.
%{\it Remarks} and Problems marked by star are not used in the sequel; although they are not necessarily complicated, they can be omitted during the first round of problem solving.
Important definitions are highlighted in \textbf{bold} for easy navigation.

\iffalse

\medskip
\textbf{On presentation at the Summer conference}
%not for arxiv? 
 
Participants are welcomed to \emph{consult} the jury on any questions on the project. 
Participants who successfully work on the project can get interesting \emph{extra problems}.

A \textit{team} working on this project may consist of any number of participants. 
For every solution \textit{written for a user} marked with either~`$+$' or `$+.$' a team gets five `beans' (see recommendations in p.~3, `How to write a proof for a user' of 
%\linebreak
\url{https://www.mccme.ru/circles/oim/multicomb.pdf}).
The jury may award extra beans for some solutions (beautiful, of hard problems, or typeset in \TeX).
The jury has infinitely many beans.
Every team initially has ten beans. 
A team may submit a solution \textit{in oral form} or as \textit{written for a developer} (i.e. less carefully than for a user) if the team has some beans. 
A team loses a bean with every attempt (successful or not).

\fi

\section{Winding number: definition and discussion}\label{s:wind}

In the plane let $O,A,B,A_1,\ldots,A_m$ be points.  

Assume that $A\ne O$ and $B\ne O$ (but possibly $A=B$). 
Recall that the \emph{oriented (a.k.a. directed) angle} $\angle AOB$ is the number $t\in(-\pi,\pi]$ such that the vector $\overrightarrow{OB}$ is codirected to the vector obtained from $\overrightarrow{OA}$ by the rotation through $t$.   
(If you are familiar with complex numbers, you can regard vectors in the plane as complex numbers, and rewrite this condition as $\overrightarrow{OB}\upuparrows e^{it}\overrightarrow{OA}$.) 

A {\bf polygonal line} $A_1\ldots A_m$ is the (ordered) set \ $(\ A_1A_2,\ A_2A_3,\ \ldots,\ A_{m-1}A_m\ )$ \ of segments. 
A {\bf closed polygonal line} $A_1\ldots A_m$ is the set \ $(\ A_1A_2,\ A_2A_3,\ \ldots,\ A_{m-1}A_m,\ A_mA_1\ )$ \ of segments.\footnote{The set of segments is not the same as the union of segments. 
Thus, strictly speaking, the polygonal line (defined here) is not a subset of the plane. 
So `oriented' or `non-oriented' is not formally applicable to polygonal lines.  
Still, we sometimes work with the set of segments as with their union, e.g. we write `a polygonal line not passing through a point'. 
The notion of polygonal line defined here is close to what is sometimes understood as `oriented polygonal line'.} 

Let $A_1\ldots A_m$ be a closed polygonal line not passing through $O$. 
The {\bf winding number} $w(A_1\ldots A_m,O)$ of $A_1\ldots A_m$ around $O$ is the number of revolutions during the rotation of vector whose origin is $O$, and whose endpoint goes along the polygonal line in positive direction.   
Rigorously, 
$$2\pi \cdot w(A_1\ldots A_m,O) := 
\angle A_1OA_2+\angle A_2OA_3+\ldots+\angle A_{m-1}OA_m+\angle A_mOA_1$$
%where angles are oriented. 
is the sum of the oriented angles.
% https://math.stackexchange.com/questions/4112869/oriented-angles-in-euclidean-geometry
%https://www.cip.ifi.lmu.de/~grinberg/NeubergMineur.pdf
%p. 4, Difference exchange formula. 
%For any four lines a; b; c; d; we have  (a; b) - (c; d) =  (a; c) - (b,d) 
%This is also called {\it the degree} of $O$ w.r.t. $A_1\ldots A_n$.
%the point w.r.t. the polygonal line.) \deg_{A_1\ldots A_n}O = A_1\ldots A_n\cdot O

\begin{figure}[!htb]
\includegraphics[scale=0.7]{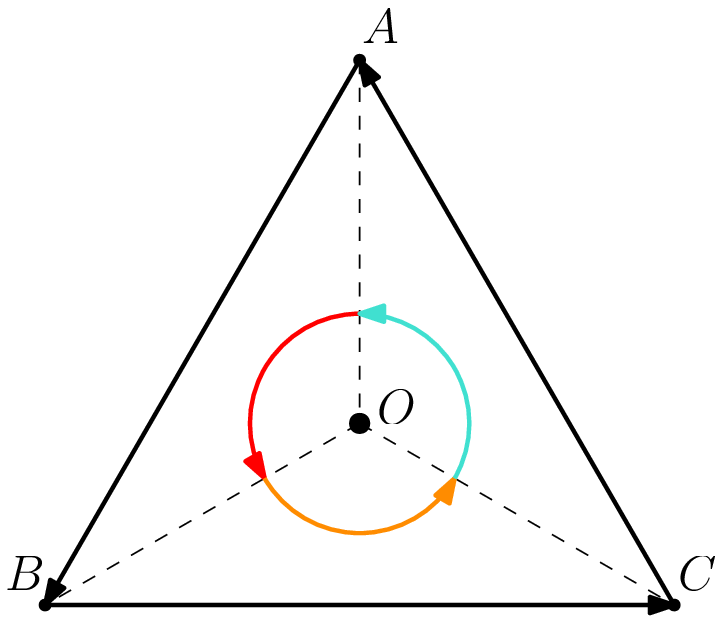}\qquad
\includegraphics[scale=0.85]{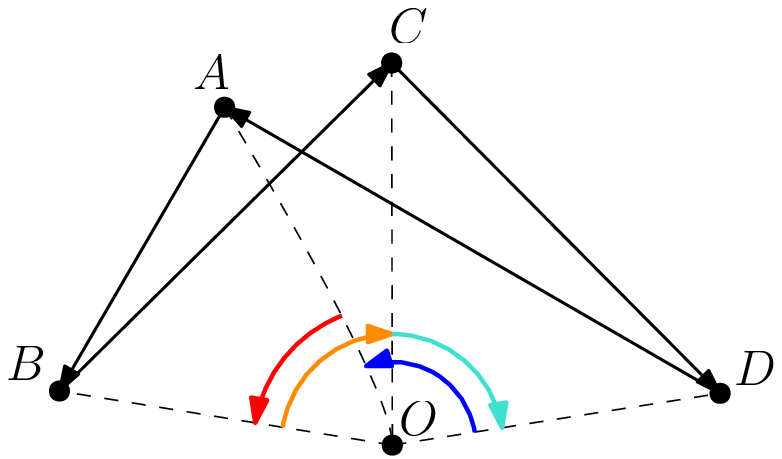}
\caption{$w(ABC, O) = +1 $ and $w(ABCD,O)=  0$}
\label{f:abco}
\end{figure}
 
E.g. in Figure \ref{f:abco}
$$
w(ABC, O) = \dfrac{1}{2\pi} \left(\angle AOB + \angle BOC + \angle COA \right) = +1\quad\text{and}
$$
$$
w(ABCD,O)=\dfrac{1}{2\pi} \left( \angle AOB + \angle BOC + \angle COD + \angle DOA \right) = \dfrac{1}{2\pi} \left( \angle BOD + \angle DOB \right) = 0.
$$

\begin{pr}\label{p:noncl} The winding number $w$ is an integer.
\end{pr}

\begin{figure}[!htb]
\includegraphics[scale=0.75]{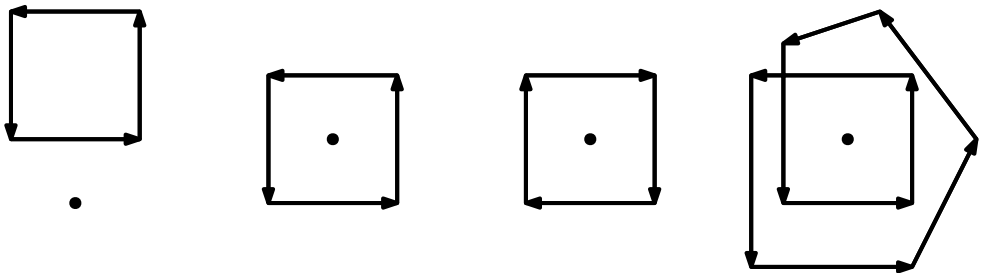}
\caption{The winding numbers equal $0,~+1,~-1,~+2$}
\label{f:wn_examples}
\end{figure}

\begin{pr}\label{p:winsim} (a) The winding number of (the outline of) any convex polygon 
%(numbered counterclockwise) 
around any point in its exterior (respectively interior) is $0$ (respectively $\pm 1$).
See Figure \ref{f:wn_examples}.  

(b) Let $ABC$ be a regular triangle and $O$ its center. 
Find $w(ABCABC,O)$. 
 
(c) For any integer $n$ and any point $O$ in the plane there is a closed polygonal line whose winding number around $O$ is $n$. 

(d) Give an example of a closed polygonal line $L$ in the plane such that $w(L,O)=0$ for any point $O\in\R^2-L$.
\end{pr}

The analogue of Assertion \ref{p:winsim}.a is correct for any closed polygonal line without self-intersections.
(Depending on the exposition, this is either a corollary of the \emph{Jordan Curve Theorem}, or a lemma in its proof.)    
The result of Problem \ref{p:winsim}.b shows that winding numbers for distinct polygonal lines with the same union of their segments can be distinct. 
 
\begin{figure}[!htb]
\includegraphics[scale=0.65]{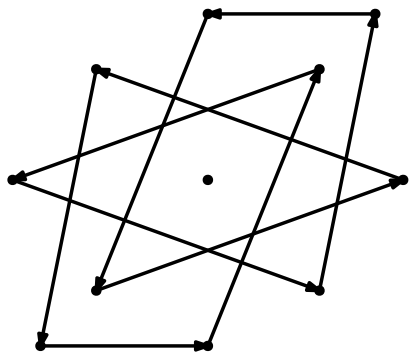}
\caption{A closed polygonal line symmetric w.r.t. a point; the winding number equals $3$}
\end{figure}
 
\begin{theorem}[Borsuk-Ulam]\label{borul-pl} 
Suppose that a closed polygonal line $A_1\ldots A_{2k}$ does not pass through a point $O$, and is symmetric w.r.t. $O$ (i.e. $O$ is the midpoint of the segment $A_jA_{k+j}$ for every $j=1,\ldots,k$).
Then the winding number of the polygonal line around $O$ is odd.
\end{theorem}

The following notion and results would be helpful. 

Let $A_1\ldots A_m$ be a polygonal line not passing through a point $O$. 
Define the real number $w'(A_1\ldots A_m,O)$ by 
$$2\pi\cdot w'(A_1\ldots A_m,O) := \angle A_1OA_2+\angle A_2OA_3+\ldots+\angle A_{m-1}OA_m.$$
Clearly, 

$\bullet$ $2\pi w(A_1\ldots A_m,O) = 2\pi w'(A_1\ldots A_m,O)+\angle A_mOA_1$; 

$\bullet$ if points $A_2,\ldots,A_{m-1}$ lie in the interior of the angle $\angle A_1OA_m$, then $w'(A_1\ldots A_m,O) = \angle A_1 O A_m$. 

\begin{pr}\label{p:w'}
(a) We have $\angle A_1OA_m = 2\pi w'(A_1\ldots A_m,O)+2\pi k$ for some integer $k$.

(b) We have $w(A_1\ldots A_m, O) = w'(A_1\ldots A_j, O) + w'(A_j\ldots A_mA_1, O)$ for every $j = 1, \ldots, m$.
\end{pr}

Denote by $\overline{l}$ the polygonal line obtained from a polygonal line $l$ by passing in the reverse direction. 

\begin{pr}\label{k23} In the plane let $O$, $A$, $B$ be three pairwise distinct points. 

(a) Let $l_1$, $l_2$, $l_3$ polygonal lines joining $A$ to $B$, and not passing through $O$.
Then 
$$w(l_1\overline{l_2}, O) + w\left(l_2\overline{l_3}, O\right) = w(l_1\overline{l_3}, O).$$

(b) For any three integers $n_1, n_2, n_3$ such that $n_1 + n_2 = n_3$ there are three polygonal lines $l_1, l_2, l_3$ joining $A$ to $B$, not passing through $O$, and such that
$$w(l_1\overline{l_2}, O) = n_1, \qquad w\left(l_2\overline{l_3}, O\right) = n_2 \quad\text{and}\quad w(l_1\overline{l_3}, O) = n_3.$$ 
\end{pr}

In Assertion \ref{p:rel}.a, Problems \ref{k23}.b, \ref{p:rel}.bc, and other \emph{examples} in this text (as opposed to \emph{assertions} distinct from \ref{p:rel}.a) you may give an heuristic rather than rigorous proof, unless you or 
%the jury 
your advisor realize that this leads to a confusion.  

\begin{pr}\label{p:3rays} Let $A_1A_2A_3$ be a regular triangle, and $O$ its center.
For $m=1,2,3$ let $l_m$ be a polygonal line disjoint with the ray $OA_m$, and joining $A_{m+1}$ to $A_{m+2}$, where the numbering is modulo 3. 
Then $w(l_1l_2l_3,O) = \pm 1$. 
\end{pr}

%\newpage
%%%%%%%%%%%%%%%%%%%%%
\section{Winding number and intersections}\label{s:wint}

This section is only used in sketches of proofs of Theorems \ref{t:radonae} and \ref{t:k5-e}.a in \S\ref{s:alem}.

\begin{pr}\label{p:rel} 
(a) In the plane points $P_0,P_1$ are joined by a polygonal line disjoint from a closed polygonal line $L$. 
Then  $w(L,P_0)=w(L,P_1)$. 

\emph{Hint:} use \emph{considerations of continuity}. 

\begin{figure}[ht]\centering
\includegraphics[scale=.8]{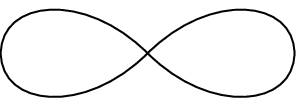} \qquad \includegraphics[scale=.8]{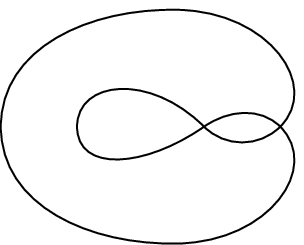} \qquad \includegraphics[scale=.8]{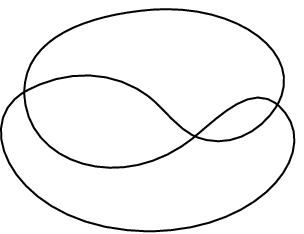} \qquad \includegraphics[scale=.8]{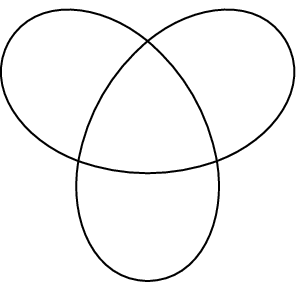} \qquad \includegraphics[scale=.8]{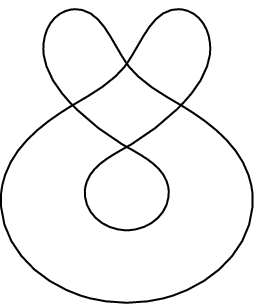}
\caption{Some closed polygonal lines}\label{f:tohu}
\end{figure}

(b) Take closed polygonal lines $L$ in the plane shown in Figure~\ref{f:tohu} (with some orientations; any of the polygonal lines does not go twice through any segment, and does not significantly change  its direction at any point).
Color the complement $\R^2-L$ by parity of the winding number of $L$. 

(c) For closed polygonal lines $L$ from (b) color the complement $\R^2-L$ by the winding number of $L$. 
\end{pr}

\begin{pr}\label{p:stokes} Take a closed and a non-closed polygonal lines $L$ and $P$ in the plane, all whose vertices are pairwise distinct and are in general position.
Let $P_0$ and $P_1$ be the starting point and the endpoint of $P$.
Assume that $P_0,P_1\not\in L$. 

(a) Then $|L\cap P|\equiv w(P_1,L)-w(P_0,L)\mod2$. 
%In other words, the ends of a polygonal line $P$ whose vertices together with the vertices of $L$ are in general position have the same parity of winding numbers if and only if $|L\cap P|$ is even.
(This is a discrete version of the \emph{Stokes theorem}.) 
 
\emph{Hint.} It suffices to prove this fact for a sufficiently small segment $P'$ such that $P'\subset P$. 

%$\bullet$ the length of $P'$ is less than some positive real number $\varepsilon(L, P)$.
 
(b) If $P_1$ is sufficiently far from $L$, then $w(P_0,L)\equiv|L\cap P|\mod2$, and $w(P_0,L)$ equals the sum of signs of intersection points of $P$ and $L$ (defined in \cite[\S1.3]{Sk18}). 
\end{pr}

The modulo 2 version of (b) follows from (a).
%Consider a ray $OB$ in general position with the line (i.e. the ray does not contain any vertex and any self-intersection point of the line). Then the winding number $w(O, A_1,\ldots A_m)$ equals the sum of signs of intersection points of the ray and the line.

Let $L$ be a closed polygonal line in the plane, all whose vertices are pairwise distinct and are in general position.  
By Assertion \ref{p:stokes}.a the complement to $L$ has a `chess-board' coloring, i.e. a coloring such that the adjacent domains have different colors. 
The \emph{modulo two interior} of $L$ is the union of black domains for a chess-board coloring (provided the infinite domain is white).
In other words, this is the set of all points $x\in\R^2-L$ for which there is a polygonal line $P$ 

$\bullet$ joining $x$ to a point `far away' from $L$ (i.e. outside the convex hill of $L$), 

$\bullet$ intersecting $L$ at an odd number of points, and %!ЕБ: an->the

$\bullet$ such that all the vertices of $L$ and $P$ are pairwise distinct, and are in general position. 

This is well-defined by \cite[Parity Lemma 1.3.3]{Sk18}.

%%%%%%%%%%%%%%%%%%%%%
\section{Winding numbers of graph drawings}\label{s:wingra}

{\bf Remark} (some rigorous definitions). 
You can work with the notions defined here at an intuitive level, before you or 
%the jury 
your advisor realize that this leads to a confusion.  

A (finite) {\bf graph} $(V,E)$ is a finite set $V$ together with a collection $E\subset {V\choose 2}$
of two-element subsets of $V$ (i.e., of non-ordered pairs of distinct elements). (The common term for this notion is {\it a graph without loops and multiple edges} or {\it a simple graph}.)
The elements of this finite set $V$ are called {\it vertices}.
%Unless otherwise indicated, we assume that $V=[|V|]$.
The pairs of vertices from $E$ are called {\it edges}.
The edge joining vertices $i$ and $j$ is denoted by $ij$ (not by $(i,j)$ to avoid confusion with ordered pairs).
A cycle in a graph is denoted by listing its vertices in their order (without commas). 
 
Informally speaking, a graph is planar if it can be drawn `without self-intersections' in the plane.
Rigorously, a graph $(V,E)$ is called {\bf planar} (or piecewise-linearly embeddable into the plane) if in the plane there exist

$\bullet$ a set of $|V|$ points corresponding to the vertices, and

$\bullet$ a set of non-self-intersecting polygonal lines joining pairs (of points) from $E$

such that none of the polygonal lines intersects the interior of any other polygonal line.\footnote{Then any two of the polygonal lines either are disjoint or intersect by a common end vertex.
We do not require that `no isolated vertex lies on any of the polygonal lines' because this property can always be achieved.}
 
\bigskip
Denote by

$\bullet$ $[n]$ the set $\{1, 2, \ldots, n\}$;

$\bullet$ $K_n$ the complete graph with the vertex set $[n]$;   

$\bullet$ $K_{m,n}$ the complete bipartite graph with parts $[m]$ and $[n]'$ (we denote by $A'$ a copy of $A$).  

\begin{figure}[h]\centering
\includegraphics{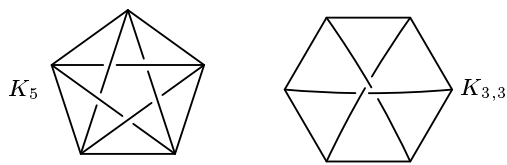} 
\caption{Non-planar graphs $K_5$ and $K_{3,3}$}\label{planar}
\end{figure}

We consider graph drawings in the plane such that the edges are drawn as polygonal lines, and intersections are allowed.
Let us give rigorous definitions.
Let $K$ be a graph with $V$ vertices. 
A (piecewise-linear) {\bf map} $f:K \to \R^2$ of $K$ to the plane is 

$\bullet$ a collection of $V$ points in the plane corresponding to the vertices, and  

$\bullet$ a collection of (non-closed) polygonal lines in the plane joining those pairs of points from the collection which correspond to pairs of adjacent vertices.\footnote{This is a slight abuse of terminology. 
A polygonal line has a starting point and an endpoint, so this is a definition of a map of an oriented graph. 
Two maps from oriented graphs are \emph{equivalent} if one of them is obtained from the other by change of orientations of some edges, and one of the corresponding collections of polygonal lines is obtained from the other by passing the corresponding polygonal lines in the reverse order.  
Rigorously speaking, a map of a graph is an equivalence class under such an equivalence relation.}

The {\bf restriction} $f|_\sigma$ to an edge $\sigma$ is the corresponding polygonal line.
The {\bf image} $f(\sigma)$ of edge $\sigma$ is the union of edges of $f|_\sigma$. 
The {\bf image} of a collection of edges is the union of images of all the edges from the collection.

\begin{figure}[!htb]
\includegraphics[scale=0.7]{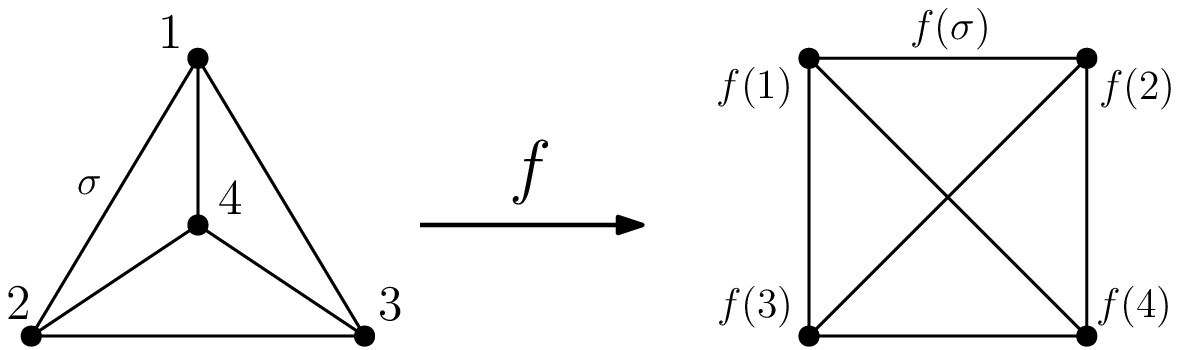}\qquad\qquad
\includegraphics[scale=0.6]{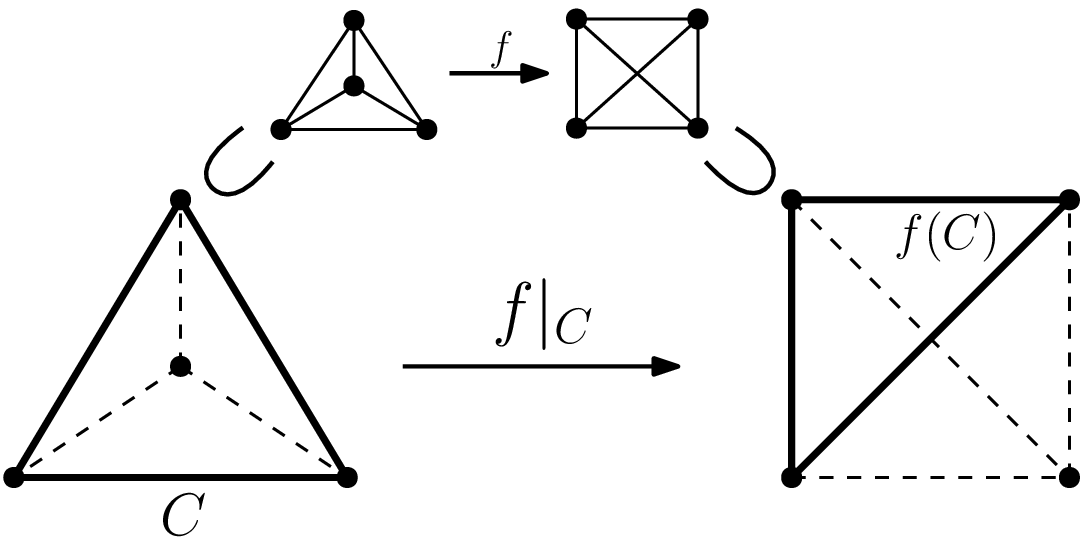}
\caption{A map $f: K_4\to\R^2$ (left) and its restriction $f|_C$ (right)} 
%where $C = 123$
\end{figure}

Let $C=v_1\ldots v_n$ be a directed (i.e. oriented) cycle in $K$. 
E.g. for $j=1,2,3,4$ denote by $C_j$ the directed cycle in $K_4$ obtained by deleting $j$ from $1234$.  
Let $f:K\to\R^2$ be a map.  
The {\bf restriction} $f|_C:C\to\R^2$ of $f$
%to $C$ 
is the closed polygonal line `formed' by the polygonal lines $f|_{v_1v_2},\ldots,f|_{v_{n-1}v_n},f|_{v_nv_1}$ in this order.    
 
\begin{pr}\label{p:rel1} For any map $f:K_4\to\R^2$ and point $O\in\R^2-f(K_4)$ we have $\sum_{j=1}^4(-1)^jw(f|_{C_j},O)=0$: 
$$-w(f|_{234},O)+w(f|_{134},O)-w(f|_{124},O)+w(f|_{123},O)=0.$$
\end{pr}

For a vertex $v$ in $K$ such that $f(v)\not\in f(C)$ denote  
$$w_f(C,v):=w(f|_C,f(v)).$$ 
 
\begin{pr}\label{p:maps} 
(a) For any integer $n$ there is a map $f:K_4\to\R^2$ such that 
$$f(1)\not\in f(C_1),\quad w_f(C_1,1)=n,\quad\text{and}\quad f(j)\not\in f(C_j),\quad w_f(C_j,j)=0\quad\text{for every}\quad j=2,3,4.$$ 

(b) For any integers $n_1,n_2,n_3,n_4$ there is a map $f:K_4\to\R^2$ such that $f(j)\not\in f(C_j)$ and $w_f(C_j,j)=n_j$ for every $j=1,2,3,4$.
\end{pr}

\begin{figure}[ht]
\centerline{\includegraphics[scale=.9]{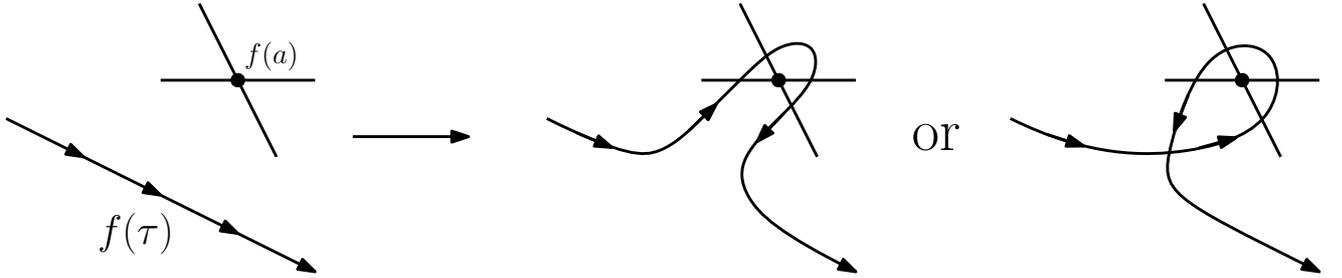}}
\caption{Finger moves (for a map $f$, of an edge $\tau$, w.r.t. the vertex $a$) of the first and the second types respectively}
\label{f:fingermove}
\end{figure}

\emph{Hint.} Transformation of a map shown in Figure \ref{f:fingermove} is useful to construct  examples. 

For more on winding number and related notions see \cite{Wn, Va81, To84, Ta88}, \cite[Theorem 2]{KK18}, \cite[\S2]{Sk18}. 

%статья из Квантика про обход фигуры без самопересечений: https://old.kvantik.com/art/files/pdf/2020-06.7-11.pdf
% XV олимпиада, заочный тур, задача 23. https://geometry.ru/olimp/2019/zaoch.pdf 

%%%%%%%%%%%%%%%%%%%%%
\section{Almost embeddings: definition and discussion}\label{s:alemdef}

\begin{theorem}[Hanani-Tutte; van Kampen]\label{t:hatuvk}
For any map $K_5\to\R^2$ there are two non-adjacent edges whose images intersect.
\end{theorem}

This follows by Assertion \ref{p:rel}.a and Theorem \ref{t:k5-e}.a below.  
(The standard proof \cite[\S1.4]{Sk18} does not use the winding number.) 
%, but is essentially the same.) , see Assertion \ref{c:ravk}.a.)
The analogue for $K_{3,3}$ holds by Assertions \ref{p:rel}.a and \ref{p:k33}.a below. 

    \begin{figure}[ht]\centering
    \includegraphics[scale=0.8]{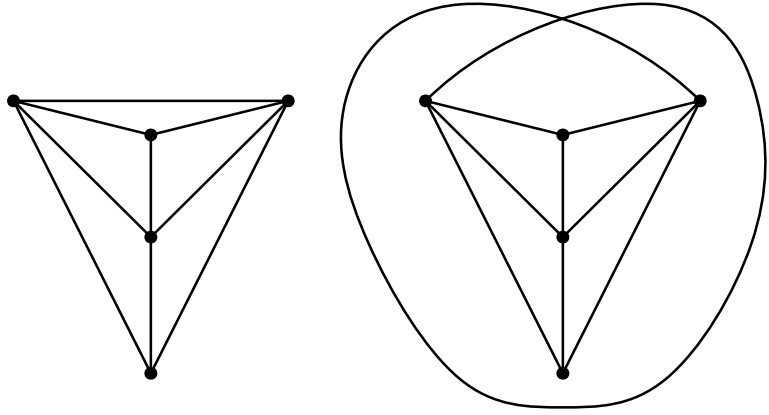}
\caption{An embedding and an almost embedding of $K_5$ without an edge}\label{k5}
\end{figure}    
 
A map $f:K\to\R^2$ of a graph $K$ is called an \textbf{almost embedding} if $f(\alpha)\cap f(\beta) = \varnothing$ for any two non-adjacent simplices (i.e. vertices or edges) $\alpha,\beta\subset K$.
In other words, if 

(i) the images of non-adjacent edges are disjoint,  

(ii) the image of a vertex is not contained in the image of any edge non-adjacent to this vertex, 
%and  
 
(iii) the images of distinct vertices are distinct. 

\smallskip
{\bf Remark.}
(a) This text primarily concerns not the problem of \emph{existence} of an almost embedding, but the \emph{invariants} of almost embeddings.  
Thus we keep in the definition the properties (ii,iii) (which can be achieved by a small enough perturbation of a map, keeping the property (i)).  
%E.g. graphs $K_5,K_{3,3}$ are not almost embeddable into the plane (this follows by Theorem \ref{t:k5-e}.a and Assertion \ref{p:k33}.a below).  
\emph{If a graph admits an almost embedding to the plane, then the graph is planar} (a proof is non-trivial). 
%and is not a part of this project). 
%This fact is non-trivial: a simple idea of proof (figure !) does not immediately work (figure !); the analogue of this fact for $Z_2$-embeddings into the sphere with four handles (figure !) is incorrect. 

(b) Almost embeddings naturally appear in topological graph theory, in combinatorial geometry, in topological combinatorics, and in studies of embeddings (of graphs in surfaces, and of hypergraphs in higher-dimensional Euclidean space).
See more motivations in \cite[\S1, ‘Motivation and background’]{ST17}, \cite[\S6.10 `Almost embeddings, $\Z_2$- and $\Z$-embeddings']{Sk}.

%%%%%%%%%%%%%%%%%%%%%%%%%%%%%%%%%
\section{Main results: winding numbers of almost embeddings}\label{s:alem}

\begin{pr}\label{p:k3} (a) For any integer $n$ and point $O$ in the plane there is an almost embedding $f:K_3\to\R^2-O$ such that $w(f|_{123},O)=n$. 

%MA??? is 4.1.b worth it? Look 4.3 (b)

(b) For any integer $n$ there is an almost embedding $f:K_4\to\R^2$ such that $w_f(123,4)=n$. 
\end{pr}

For any embedding $f:K_3\sqcup\{4\}\to\R^2$ we have $w_f(123,4)\in\{-1,0,1\}$ (this result is close to Jordan Curve Theorem).  

Recall that 
$$\sum_{j=1}^4 w_f(C_j,j) = w_f(234,1)+w_f(134,2)+w_f(124,3)+w_f(123,4).$$
 
\begin{theorem}\label{t:radonae} For any almost embedding $f:K_4\to\R^2$ we have 
$\sum_{j=1}^4 w_f(C_j,j) \equiv1\mod2$. 
\end{theorem}

The analogue of Theorem \ref{t:radonae} for embeddings instead of almost embeddings is simple (and is close to Jordan Curve Theorem).   
Moreover, for any embedding $f:K_4\to\R^2$ three of the four numbers from Theorem \ref{t:radonae} are zeroes, and the remaining one is $\pm1$. 
The analogue of Theorem \ref{t:radonae} for maps instead of almost embeddings is incorrect by Assertion \ref{p:maps}.    
Unlike Assertion \ref{p:rel1}, Theorem \ref{t:radonae} does not come from the `relation $123+134+142+243=0$ in the graph'. 
%Proofs of Theorems \ref{t:radonae}, \ref{t:k5-e} and \ref{p:k33} are not part of this project 

A \emph{general position} map $f:K_n\to\R^2$ is defined in \cite[\S1.4]{Sk18}. 

\begin{proof}[Sketch of a proof of Theorem \ref{t:radonae}] 
For a general position map $f:K_4\to\R^2$ let the \emph{Radon number} $\rho(f)\in\Z_2$ be the  the parity of the sum of 

$\bullet$ the number of intersections points of the images of non-adjacent edges, and

$\bullet$ the number of vertices $j$ whose images belong to the interior modulo 2 of $f(C_j)$. 
%the image of the cycle formed by the three edges not containing the vertex. 

By Assertion \ref{p:stokes}.a, \emph{for every general position almost embedding $f:K_4\to\R^2$ the parity of $\sum_{j=1}^4 w_f(C_j, j)$ equals $\rho(f)$}. 
Theorem \ref{t:radonae} is deduced using approximation from this result and the following celebrated topological Radon  theorem for the plane \cite[Lemma 2.2.3]{Sk18}:  
\emph{for any general position map $f:K_4\to\R^2$ the Radon number $\rho(f)$ is odd.}
\end{proof}

\begin{example}\label{e:radonae} (a) For any integer $n$ there is an almost embedding $f:K_4\to\R^2$ such that $w_f(C_1,1)=n$,  $w_f(C_j,j)=0$ for every $j=2,3$, and $w_f(C_4,4)=n+1$. 

% For any integers $m$ and $n$ there is an almost embedding $f:K_4\to\R^2$ such that 
% $$w_f(C_1, 1)=m,\quad w_f(C_2, 2)=n,\quad w_f(C_3, 3)=n-m\quad\text{and}\quad w_f(C_4, 4)=1.$$

(b) For any integers $n_1, n_2, n_3, n_4$ such that $\sum_{j=1}^4 (-1)^j n_j = \pm 1$ there is an almost embedding $f:K_4\to\R^2$ such that $w_f(C_j, j) = n_j$ for every $j = 1,2,3,4.$ 
\end{example}

\begin{example}\label{c:radonae} (a) (E. Morozov) There is an almost embedding $f:K_4\to\R^2$ such that 
\linebreak
$\sum_{j=1}^4 (-1)^j w_f(C_j, j) \ne \pm 1$.

%(b) For any integers $k_1,k_2,k_3,k_4$ there is an almost embedding $f:K_4\to\R^2$ such that 
%\linebreak
%$\sum_{j=1}^4 k_j w_f(C_j, j) \ne \pm 1$. 

(b) (conjecture; see \cite{ALM}) For any integers $n_1,n_2,n_3,n_4$ whose sum is odd there is an almost embedding $f:K_4\to\R^2$ such that $w_f(C_j, j)=n_j$ for every $j=1,2,3,4$. 
\end{example}
  
Recall that $K-e$ is the graph obtained from a graph $K$ by deleting an edge $e$. 

\begin{theorem}\label{t:k5-e} For any almost embedding $f:K_5-45\to\R^2$ we have 

(a) $w_f(123,4)-w_f(123,5)\equiv1\mod2$; \qquad 

(b)* $w_f(123,4)-w_f(123,5)=\pm 1$. 
\end{theorem}

Cf. Theorem \ref{t:radonae} and Conjecture \ref{c:radonae}.b.  
The analogue of (b) for embeddings instead of almost embeddings is simple (and is close to Jordan Curve Theorem).    
Part (a) is not hard and is well-known, while part (b) is a recent non-trivial result of \cite{Ga23}. 
 
\begin{pr}\label{p:k5-emap} 
(a) Is the analogue of Theorem \ref{t:k5-e}.a correct for maps instead of almost embeddings?  

(b) For any integer $n$ there is an almost embedding $f:K_5-45\to\R^2$ such that $w_f(123,5)=n$. 
(See Figure \ref{k5}, right, for $n=2$.)
\end{pr}

\begin{proof}[Sketch of the proof of Theorem \ref{t:k5-e}.a] 
For a general position map $f:K_5\to\R^2$ color in red the intersections points of the images of non-adjacent edges. 
Let \emph{the van Kampen number} $v(f)\in\Z_2$ be the parity of the number of red points.   
By Assertion \ref{p:stokes}.a, \emph{for any general position map $f:K_5\to\R^2$ whose restriction to $K_5-45$ is an almost embedding the parity of $w_f(123, 5)-w_f(123, 4)$ equals $v(f)$.}   
Theorem \ref{t:k5-e}.a is deduced from this and the following celebrated van Kampen-Flores theorem for the plane \cite[Lemma 1.4.3]{Sk18}: \emph{for any general position map $f:K_5\to\R^2$ the van Kampen number $v(f)$ is odd.} 
%(Cf. \cite[Remark 4]{Ga23}.)
\footnote{Theorem \ref{t:k5-e}.b is an \emph{integer version for almost embeddings} of this theorem.   
Observe that this theorem has no \emph{integer version for maps} (this is known and is explained in \cite[Remark 4]{Ga23}).}
\end{proof}

\begin{pr}\label{p:k33} Take an edge $ab$ of $K_{3,3}$. 
Denote by $C=C_{ab}$ somehow oriented cycle $K_{3,3}-a-b$ of length 4. 
For any almost embedding $f:K_{3, 3}-ab\to\R^2$ we have 

(a) $w_f(C,a)-w_f(C,b)\equiv1\mod2$; \qquad (b) (conjecture) $w_f(C,a)-w_f(C,b)=\pm 1$.
\end{pr}

\begin{figure}[ht]\centering
\includegraphics[scale=.7]{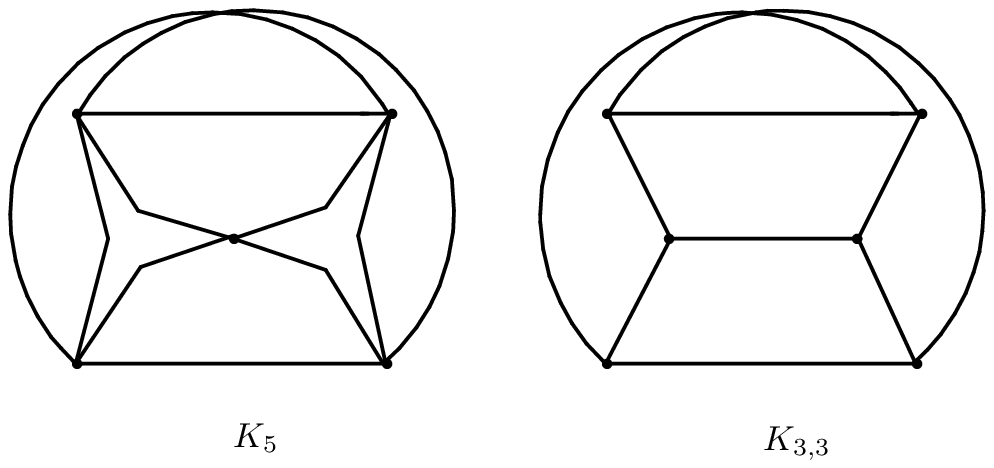}
\caption{`Almost embedding' $K_5\to K_{3,3}$}\label{fig-almemb}
\end{figure}

Part (a) is proved analogously to Theorem \ref{t:k5-e}.a. 
The analogue of (b) for embeddings instead of almost embeddings is simple (and is close to Jordan Curve Theorem).    
Beware that a direct proof of (b) might contains technical details (like for Theorem \ref{t:k5-e}.b); perhaps there is a simple a reduction to Theorem \ref{t:k5-e}.b, see Figure  \ref{fig-almemb}. 
 
\begin{pr}[open problem; riddle]\label{p:cube} Let $K$ be the graph of \qquad

(a) a cube; \qquad (b) an octahedron. 

For an almost embedding $f:K\to\R^2$ consider the collection $w_f(C,v)$ of integers, where $v\in K$ is a vertex, and $C\subset K-v$ is an oriented cycle. 
Describe collections realizable by almost embeddings $f:K\to\R^2$.
\end{pr}

%\begin{remark}\label{r:mot} %\end{remark}
\smallskip
{\bf Remark.}
(a) The integer $w_f(C,v)$ and the invariants studied in \S\ref{s:gawh} are \emph{(almost) isotopy invariants} of an (almost) embedding $f:K\to\R^2$. 
They are parts of the \emph{Haefliger-Wu invariant} of $f$ \cite[\S5]{Sk06}.  

(b) An algebraic version of almost embedings ($\Z_2$-embeddings) appeared in 1930s and is actively studied in graph theory since 2000s. 
See e.g. surveys \cite{SS13}, \cite[\S6.10 `Almost embeddings, $\Z_2$- and $\Z$-embeddings']{Sk}, and the paper \cite{FK19} relating $\Z_2$-embeddings to low rank matrix completion problem. The analogues of assertions \ref{t:radonae}, \ref{t:k5-e}.a, \ref{p:k33}.a are correct for $\Z_2$-embeddings. 
We conjecture that the analogues of assertions \ref{t:k5-e}.b, \ref{p:k33}.b are correct for $\Z$-embeddings, but are incorrect for $\Z_2$-embeddings.

(c) A \emph{hypergraph} is a higher-dimensional analogue of graph: together with edges joining pairs of points one considers triangles spanned by triples of points, etc. 
A classical problem in topology, combinatorics and computer science is to find criteria (and algorithms) for realizability of hypergraphs in Euclidean space of given dimension $d$. 

\begin{figure}[h]\centering
\includegraphics{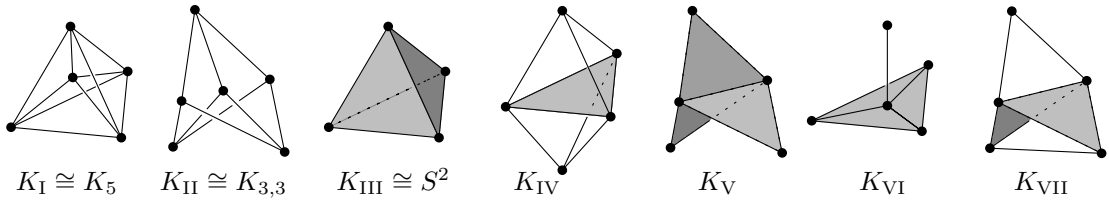}
\caption{Two-dimensional hypergraphs non-embeddable in the plane}\label{hj}
\end{figure}

Such a criterion was obtained in 1930s-1960s by classical figures in topology. 
%Egbert van Kampen, Wu Wen-tsun, Arnold Shapiro, Andr\'e Haefliger and Claude Weber. 
The criterion is stated in terms of certain configuration space, yielded many specific corollaries, and works for $2d\ge3k+3$, where $k$ is the dimension of the hypergraph. 
A polynomial algorithm based on this criterion was obtained in 2013. 
%by Martin {\v{C}}adek, Marek Kr\v{c}\'{a}l and Lucas Vok\v{r}\'{\i}nek. 
The non-existence of a polynomial algorithm for $2d<3k+2$ was announced in 2019 by 
Marek Filakovsk\'y, Ulrich Wagner and Stephan Zhechev. 
A mistake was found in 2020 by Arkadiy Skopenkov (and recognized by the authors). 
The mistake was that in a higher-dimensional analogue of Theorem \ref{t:k5-e}.b (and of Example \ref{e:k6}.a) certain linking invariant can assume values distinct from $\pm1$. 
In 2020 Roman Karasev and Arkadiy Skopenkov showed that the linking invariant can assume any odd value.
Their conjecture that the same holds for graphs in the plane was refuted by Timur Garaev, see Theorem \ref{t:k5-e}.b. 
For references see surveys \cite[\S5]{Sk06}, \cite[\S3]{Sk18}, and recent research papers \cite{KS20, Ga23}.

%\newpage
%%%%%%%%%%%%%%%
\section{Triodic and cyclic Wu numbers}\label{s:gawh}

In the plane let $l_1,l_2,l_3$ be polygonal lines joining a point $O$ to points $A_1,A_2,A_3$, respectively. 
Assume that $A_i\not\in l_j$ for every $i\neq j$.
(In other words, recall that $K_{3,1}$ is the graph with vertices $\{1,2,3,1'\}$, where $\deg 1'=3$ and $\deg m=1$ for each $m\in [3]$; take an almost embedding $f:K_{3,1}\to\R^2$ and denote $l_m:=f(1'm)$ for each $m\in [3]$.)

\begin{figure}[ht]\centering{
\includegraphics[scale=0.22]{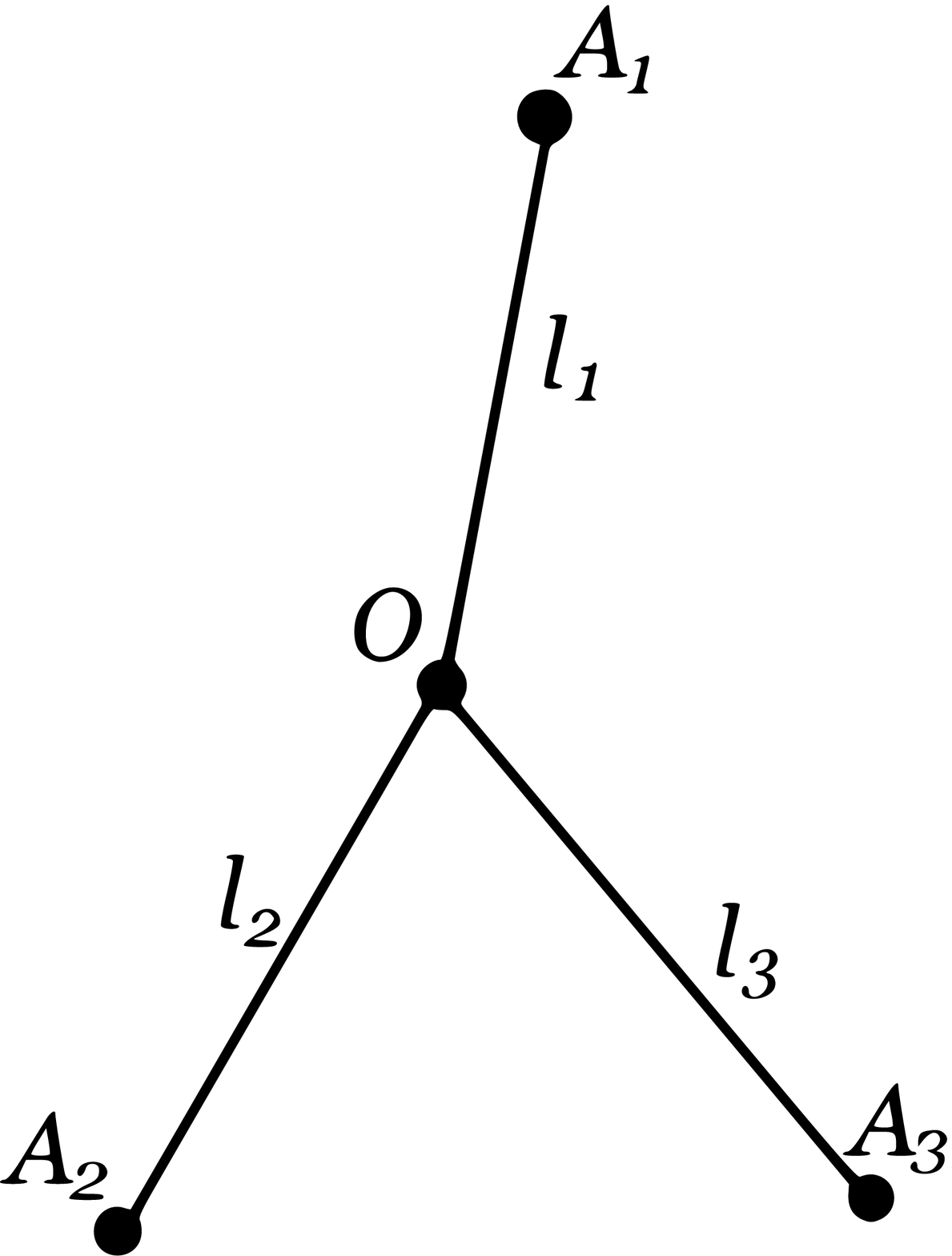} \qquad \includegraphics[scale=0.26]{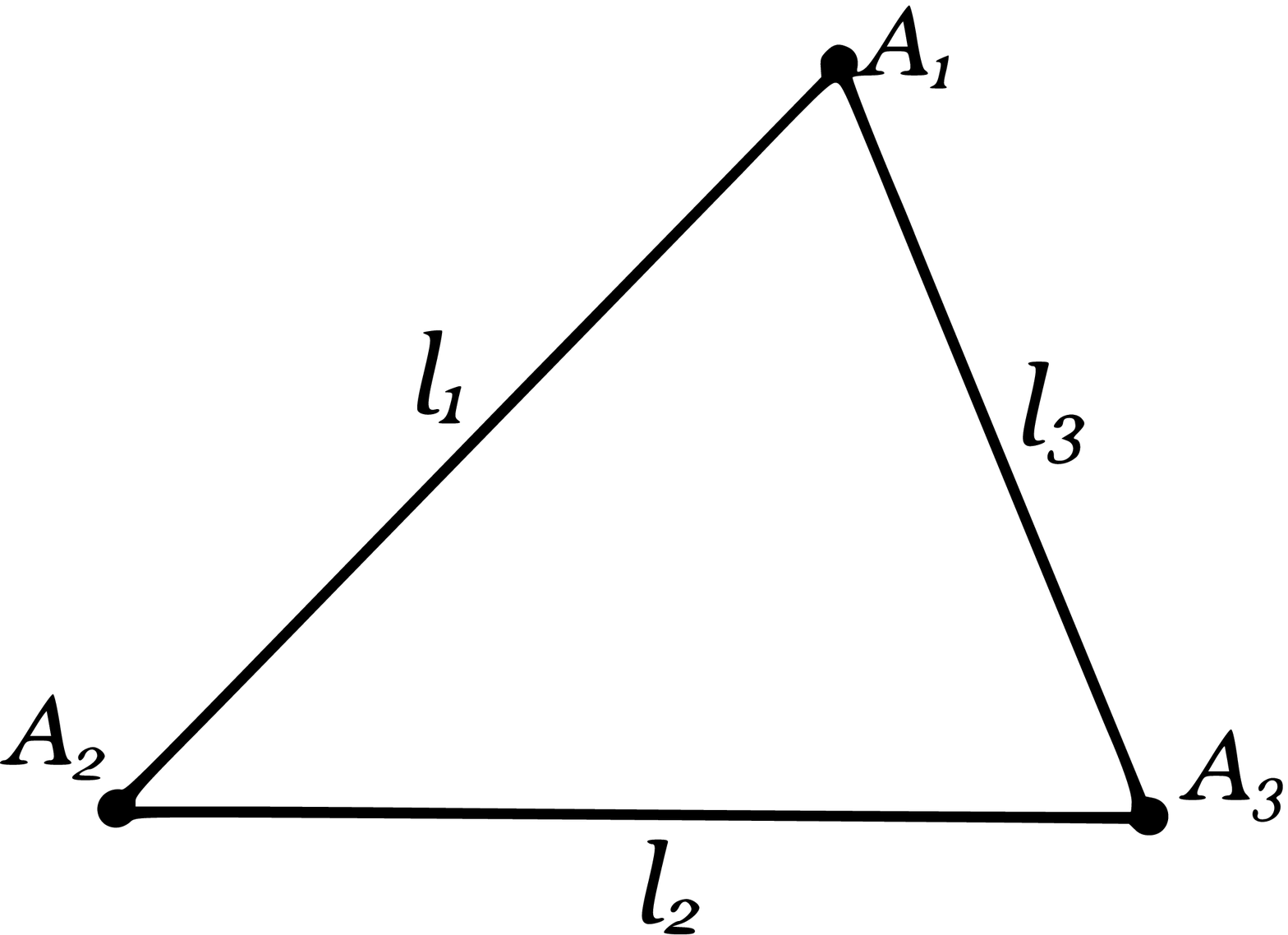}}
\caption{A triod and a triangle}
\label{f:triod}
\end{figure}

The \textit{triodic Wu number} $\wu(l_1,l_2,l_3)$ is defined to be the number of revolutions in the following rotation of vector:

%\wu_t

$\bullet$ from $\overrightarrow{A_1A_2}$ to $\overrightarrow{A_1A_3}$ as the second point of the vector moves along 
$\overline{l_2}l_3$, then

$\bullet$ from $\overrightarrow{A_1A_3}$ to $\overrightarrow{A_2A_3}$ as the first point of the vector moves along $\overline{l_1}l_2$, then

$\bullet$ from $\overrightarrow{A_2A_3}$ to $\overrightarrow{A_2A_1}$ as the second point of the vector moves along $\overline{l_3}l_1$, then

$\bullet$ from $\overrightarrow{A_2A_1}$ to $\overrightarrow{A_3A_1}$ along $\overline{l_2}l_3$, then

$\bullet$ from $\overrightarrow{A_3A_1}$ to $\overrightarrow{A_3A_2}$ along $\overline{l_1}l_2$, then

$\bullet$ from $\overrightarrow{A_3A_2}$ to $\overrightarrow{A_1A_2}$ along $\overline{l_3}l_1$.

This equals twice the (non-integer) number of revolutions in first three rotations above.
Rigorously,
$$\wu(l_1, l_2, l_3) := 
w'(\overline{l_2}l_3, A_1)+w'(\overline{l_1}l_2, A_3)+w'(\overline{l_3}l_1, A_2)+w'(\overline{l_2} l_3, A_1)+w'(\overline{l_1} l_2, A_3)+w'(\overline{l_3} l_1, A_2) =$$ 
$$=2\left(w'(\overline{l_2} l_3, A_1)+w'(\overline{l_1} l_2, A_3)+w'(\overline{l_3} l_1, A_2)\right). \qquad(*)$$

\begin{figure}[ht]\centering 
        \includegraphics[scale=0.5]{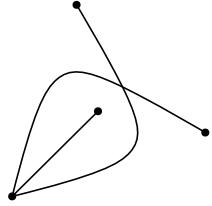}
	\caption{Three polygonal lines whose triodic Wu number equals $\pm3$}
\label{f:k31}
\end{figure}  
 
\begin{pr}\label{e:triod} (a) For the three segments joining the vertices $A_1,A_2,A_3$ of a regular triangle to its center $O$, the triodic Wu number equals $\pm1$. 

(b) For the three polygonal lines shown in Figure \ref{f:k31} the triodic Wu number equals $\pm3$. 
\end{pr}

\begin{pr}\label{p:triod} (a) For any polygonal lines $l_1,l_2,l_3$ as above (i.e. for any an almost embedding $f:K_{3,1}\to\R^2$) the triodic Wu number is odd. 

%move to {e:triod}.c?
(b) For any integer $n$ there are polygonal lines $l_1,l_2,l_3$ as above (i.e. there is an almost embedding $f:K_{3, 1}\to\R^2$) whose triodic Wu number equals $2n+1$.
%(See e.g. Figure \ref{f:k31}.) 

(c) For any embedding $f:K_{3,1}\to\R^2$ the triodic Wu number is $\pm1$. 

(d) Permutation of the polygonal lines $l_1,l_2,l_3$ as above multiplies the triodic Wu number by the sign of the permutation.  
\end{pr}

In the plane let $A_1,A_2,A_3$ be points and $l_1,l_2,l_3$ polygonal lines joining $A_1$ to $A_2$, $A_2$ to $A_3$, $A_3$ to $A_1$, respectively (and thus forming a closed polygonal line).  
Assume that $A_i$ is not contained in $l_{i+1}$ for each $i=1,2,3$ (the numeration is modulo 3; in other words, the polygonal lines form an almost embedding $K_3\to\R^2$). 
The \textit{cyclic Wu number} $\wu(l_1,l_2,l_3)$ is defined to be twice the number of revolutions in the following rotation of vector:  

$\bullet$ from $\overrightarrow{A_1A_2}$ to $\overrightarrow{A_1A_3}$, as the second point of the vector moves along $l_2$, then

$\bullet$ from $\overrightarrow{A_1A_3}$ to $\overrightarrow{A_2A_3}$, as the first point of the vector moves along $l_1$, then

$\bullet$ from $\overrightarrow{A_2A_3}$ to $\overrightarrow{A_2A_1}$, as the second point of the vector moves along $l_3$. 

In other words, $\wu(l_1, l_2, l_3)$ is defined by the following formula analogous to (*):   
$$\wu(l_1,l_2,l_3) := 2\left(w'(l_2, A_1)+w'(l_1, A_3)+w'(l_3, A_2)\right). \qquad(**)$$
  
\begin{pr}\label{p:off} (a') If the three polygonal lines $l_1,l_2,l_3$ as above are sides of a triangle, then the cyclic Wu number is $\pm1$. 

(a) For any polygonal lines $l_1,l_2,l_3$ as above the cyclic Wu number is odd. 

(b) For any integer $n$ there are polygonal lines $l_1,l_2,l_3$ as above whose  cyclic Wu number equals $2n+1$.
\end{pr}

If three polygonal lines as above form a simple closed polygonal line, then the cyclic Wu number is $\pm1$ (this is close to Jordan Curve Theorem).    

The cyclic Wu number is similar to, but is distinct from, the \emph{degree} of a closed curve. 
%\cite{}. 

\begin{pr}\label{p:wu-conj} (a) For any almost embedding $f:K_4\to\R^2$ we have 
$$\wu(f|_{12}, f|_{23}, f|_{31}) + \wu(f|_{41}, f|_{42}, f|_{43}) = 2w_f(123,4).$$ 

(b) (conjecture; see \cite{Za}) For any two odd integers $n, m$ there is an almost embedding $f:K_4 \to \R^2$ such that $\wu(f|_{12}, f|_{23}, f|_{31}) = m$ and $\wu(f|_{41}, f|_{42}, f|_{43}) = n$. 
\end{pr}

\begin{pr}\label{p:wi-wu} (a) For any almost embedding $f:K_5-45\to\R^2$ we have 
$$\wu(f|_{41}, f|_{42}, f|_{43}) - \wu(f|_{51}, f|_{52}, f|_{53}) = 2(w_f(123,4)-w_f(123,5)).$$ 
 
(b) (conjecture) For any almost embedding $f:K_{3,2}\to\R^2$ we have  
$$\wu(f|_{1'1}, f|_{1'2}, f|_{1'3}) - \wu(f|_{2'1}, f|_{2'2}, f|_{2'3}) = \pm2.$$ 
\end{pr}

By (a), a simple proof of (b) would give a simple proof of Theorem \ref{t:k5-e}.b.

%%%%%%%%%%%%%%%%%%%%%%
\section{3-dimensional analogues}\label{s:3dim}

%Let us present 3-dimensional analogues of the results from \S\ref{s:alem}. 
%\begin{theorem}[
\emph{The Linear Conway--Gordon--Sachs Theorem.} 
%asserts that 
%]\label{rlt-cgslin}
If no 4 of 6 points in 3-space lie in one plane, then there are two linked triangles with vertices at these 6 points (i.e. the first triangle intersects the outline of the second triangle exactly at one point).
%\end{theorem}
For a proof see survey \cite{Sk14}. 
%; this is not a part of this project. 

\begin{figure}[ht]\centering
\includegraphics{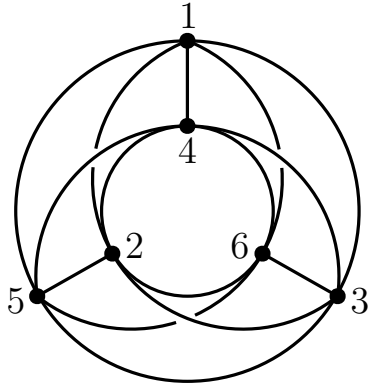}
\caption{A projection on the plane of an embedding $K_6 \to \R^3$}
\label{f:k6-emb}
\end{figure}

Equivalent rigorous definitions of a \emph{linking number} of disjoint closed polygonal lines in 3-space  
%(used in the following problem) 
can be found in \cite[\S\S 4,8]{Sk20u}, \cite[\S\S 1.2,1.3]{Sk24}, \cite[\S\S 4.2,4.3]{Sk}. 
%They will be discussed with participants successfully solving other problems.  

\begin{example}[cf. Figure \ref{f:k6-emb}]\label{e:k6}* 
(a) For any integer $n$ there are six points in 3-space, and non-self-intersecting polygonal lines joining each pair of them, and such that 

\quad (a1) the interior of one polygonal line is disjoint with any other polygonal line, 

\quad (a2) the linking number of one (unordered) pair of disjoint cycles of length 3 formed by the polygonal lines is $2n+1$, and 

\quad (a3) the linking number of any other pair of disjoint cycles of length 3 formed by the polygonal lines is zero.  
 
(b) (conjecture) Take any 10 integers $n_{123,456},n_{124,356},\ldots$ corresponding to the 10 non-ordered splittings of $[6]$ into two 3-element subsets. 
If the sum of the integers is odd, then there are 6 points $1,2,3,4,5,6$ in 3-space, and non-self-intersecting polygonal lines joining each pair them, for which (a1) holds, and the linking number of every pair $\{ijk,pqr\}$ of disjoint cycles  of length 3 formed by the polygonal lines equals $n_{ijk,pqr}$.
\end{example}

%https://arxiv.org/pdf/2105.10984, \S2.3, On Figure 3(c). 

Part (a) is proved in \cite[Proposition 1.2]{KS20} but could have been known before. 
Part (b) could perhaps be proved using \cite[Lemma 2.1]{KS20}. 
For higher-dimensional analogues see survey \cite{Sk14} and \cite[\S1]{KS20}. 
In those results the linking number can assume any odd value, like Conjecture \ref{c:radonae}.b, but unlike Conjecture \ref{p:k33}.b and Theorem \ref{t:k5-e}.b.

%%%%%%%%%%%%%%%%%%%%%%%
%\newpage
\section{Answers, hints and solutions}\label{s:answ} 
%presented at the intermediate discussion

\textbf{\ref{p:noncl}.}
For $j=1,\ldots,m$ let $t_j:=\angle A_jOA_{j+1}$, where $A_{m+1}=A_1$.  
Then  
$$\overrightarrow{OA_1} \upuparrows e^{it_m}\overrightarrow{OA_m} \upuparrows e^{i(t_m+t_{m-1})}\overrightarrow{OA_{m-1}} \upuparrows \ldots \upuparrows e^{i(t_m+t_{m-1}+\cdots+t_1)}\overrightarrow{OA_1}.$$
Hence $(t_m+t_{m-1}+\cdots+t_1)/2\pi$ is an integer. 

\smallskip
\textbf{\ref{p:winsim}.} (a) Let $\Omega$ be a convex polygon. 

If $O$ is in the exterior of $\Omega$, draw two supporting 
%https://en.wikipedia.org/wiki/Supporting_line
lines from $O$ to $\Omega$. 
Take two points $A,B$ from intersections of the lines and $\partial \Omega$. 
Then $w(\partial \Omega,O) = \frac{1}{2\pi} (\angle AOB + \angle BOA) = 0$. 

If $O$ is in the interior of $\Omega$, draw a regular triangle $ABC$ centered at $O$.
Take three intersection points $A'$, $B'$, $C'$ of the rays $OA$, $OB$, $OC$ with $\partial \Omega$.
They split $\partial \Omega$ in three polygonal lines.
We have
$$w(\partial \Omega,O) = \frac{1}{2\pi} (\angle A'OB' + \angle B'OC' + \angle C'OA') = \frac{3}{2\pi} \angle A'OB' = \frac{3}{2\pi} \angle AOB = \pm 1.$$
 
(b) By (a), we have $w(ABCABC,O) = 2 \cdot w(ABC,O)= \pm 2$.
  
(c) If $n=0$, example is a single-point closed polygonal line. 
If $n\ne 0$, let $ABC$ be a regular triangle oriented clockwise if $n<0$, and counterclockwise in the opposite case. 
Let $O$ be its center. 
Consider the closed polygonal line $L = \underbrace{ABC\ldots ABC}_{|n| \text{ times}}$.
By a refinement of (a), we have $w(L,O) = |n| \cdot w(ABC,O) = n$. 

(d) A trivial example is a single-point closed polygonal line. 
Another example: a closed polygonal line $ABCB$ for any three points $A, B$ and $C$ in the plane.       

\smallskip
\textbf{\ref{borul-pl}.}
By the symmetry, $w'(A_1 \ldots A_{k+1},O) = w'( A_{k+1} \ldots A_{2k}A_1,O)$.
Then 
$$
w(A_1\ldots A_{2k},O) = w'(A_1 \ldots A_{k+1},O) + w'( A_{k+1} \ldots A_{2k}A_1,O) =
$$
$$
= 2w'(A_1 \ldots A_{k+1},O) = 2 \left( \frac{\angle A_1OA_{k+1}}{2\pi} + n \right) = 1 + 2n
$$
for some integer $n$. 
Here the last but one equality follows from Assertion~\ref{p:w'}.a.

\smallskip
\textbf{\ref{p:w'}.}
(a) Using the formula from the proof of Assertion~\ref{p:noncl}, we have 
$$
e^{i\angle A_1OA_m}\overrightarrow{OA_1}  \upuparrows \overrightarrow{OA_m} \upuparrows e^{2\pi 
 i w'(A_1 \ldots A_m, O)}\overrightarrow{OA_1}.
$$
Hence $\angle A_1OA_m - 2\pi  w'(A_1 \ldots A_m, O) = 2\pi k$ for some integer $k$.

\smallskip
\textbf{\ref{k23}.} \emph{Hint.} (a) Use Assertion \ref{p:w'}.b and the following equality: $w'(l, O) = -w'(\overline{l}, O)$ for every point $O$ and every polygonal line $l$ not passing through $O$.

\smallskip
\textbf{\ref{p:3rays}.} \emph{Hint.}
We can assume that the vertices of the triangle are numbered counter-clockwise.
%, and for any $m \in [3]$ the polygonal line $l_m$ joins the point $A_{m+1}$ to $A_{m+2}$, where $A_4 = A_1$. 
Prove that $w'(l_0,O) = w'(l_1,O) = w'(l_2,O) = \frac{2\pi}{3}$.
For that denote by $B_1 \ldots B_m$ the sequence of vertices of $l_0$. 
For $j \in [m-1]$ let $t_j := \angle B_jOB_{j+1} \in (-\pi, \pi)$.
Prove that for every $j \in [m-1]$ we have $T_j := t_1 + \ldots + t_j \in (-\frac{2\pi}{3}, \frac{4\pi}{3})$, and $w'(l_0,O) = T_{m-1} = \frac{2\pi}{3} + 2\pi k$ for some integer $k$. 
Then deduce that $k = 0$. 

% Denote by $l_m^\perp$ a line perpendicular to $l_m$ and passing through $O$.
% For $m \in [3]$ use \ref{slight}.b to move every point of $l_m$ to the half-plane of $l_m^\perp$ opposite to the ray $OA_m$.
% Then use \ref{slight}.b to make $l_m$ into a segment.
% Thus, $w(O,l_1l_2l_3) = w(O, A_1A_2A_3) = \pm 1$.

\begin{figure}[ht]\centering
\includegraphics[scale=.8]{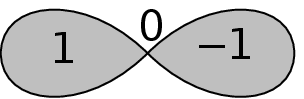} \qquad \includegraphics[scale=.8]{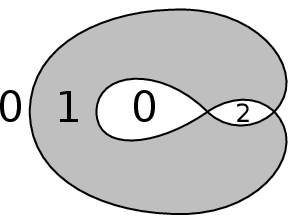} \qquad \includegraphics[scale=.8]{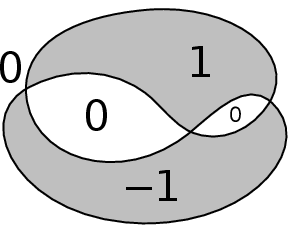} \qquad \includegraphics[scale=.8]{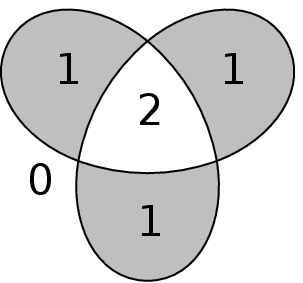} \qquad \includegraphics[scale=.8]{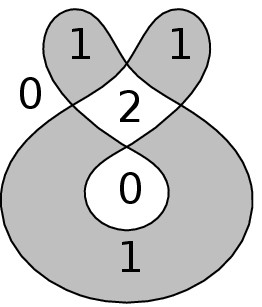}
\caption{Coloring of the complement of the closed polygonal lines from Figure~\ref{f:tohu} by the winding number}
\label{f:tohu_colored}
\end{figure}

\smallskip
\textbf{\ref{p:rel}.} (b,c) Figure \ref{f:tohu_colored}.

\smallskip
\textbf{\ref{p:rel1}.} \emph{Hint.}
Use the following equalities: $w'(f|_{ij},O) =- w'( f|_{ji},O)$ for every edge $ij$, and 
$w(f|_{ijk},O)=w'(f|_{ij},O)+w'(f|_{jk},O)+w'(f|_{ki},O)$ for every cycle $ijk$. 
%(Here we use oriented edges)
%$$w(O, f|_C) = \sum\limits_{ij \text{ is an edge of $C$}} w'(O, f|_{ij})$$ for every cycle $C$.

\smallskip
\textbf{\ref{p:maps}.} \emph{Hints.} (a)  A square with the diagonals forms a map $f:K_4\to\R^2$; assume the vertices $f(i),~ i \in [4]$ are numbered counterclockwise.
Then $w_f(C_j,j) = 0$ for every $j \in [4]$.
%It is also clear that winding number does not change under subdivision of any edge.
Make $|n|$ finger moves (Figure~\ref{f:fingermove}) of the edge $24$ w.r.t. the vertex $1$ of the first/second type if $n$ is positive/negative respectively.
The obtained map $f_1$ is as required.

(b) Consider the map $f_1$ from the proof of (a) for $n = n_1$.
Make $|n_2|$ finger moves of the edge $13$ w.r.t. the vertex $2$ of the first/second type if $n_2$ is positive/negative respectively.
Denote the obtained map by $f_2$.
Make $|n_3|$ finger moves of the edge $24$ w.r.t. the vertex $3$ of the first/second type if $n_3$ is positive/negative respectively.
Denote the obtained map by $f_3$.
Make $|n_4|$ finger moves of the edge $13$ w.r.t. the vertex $4$ of the first/second type if $n_4$ is positive/negative respectively.
The obtained map $f_4$ is as required.

%(cf. Figure \ref{f:maps}).
%\begin{figure}[ht]\centering
%\includegraphics[scale=.2]{sol_emb_k4.eps}  
%\caption{A map $f:K_4\to\R^2$ with given winding numbers}\label{f:maps}
%\end{figure}

\begin{figure}[ht]\centering
\includegraphics[scale=.16]{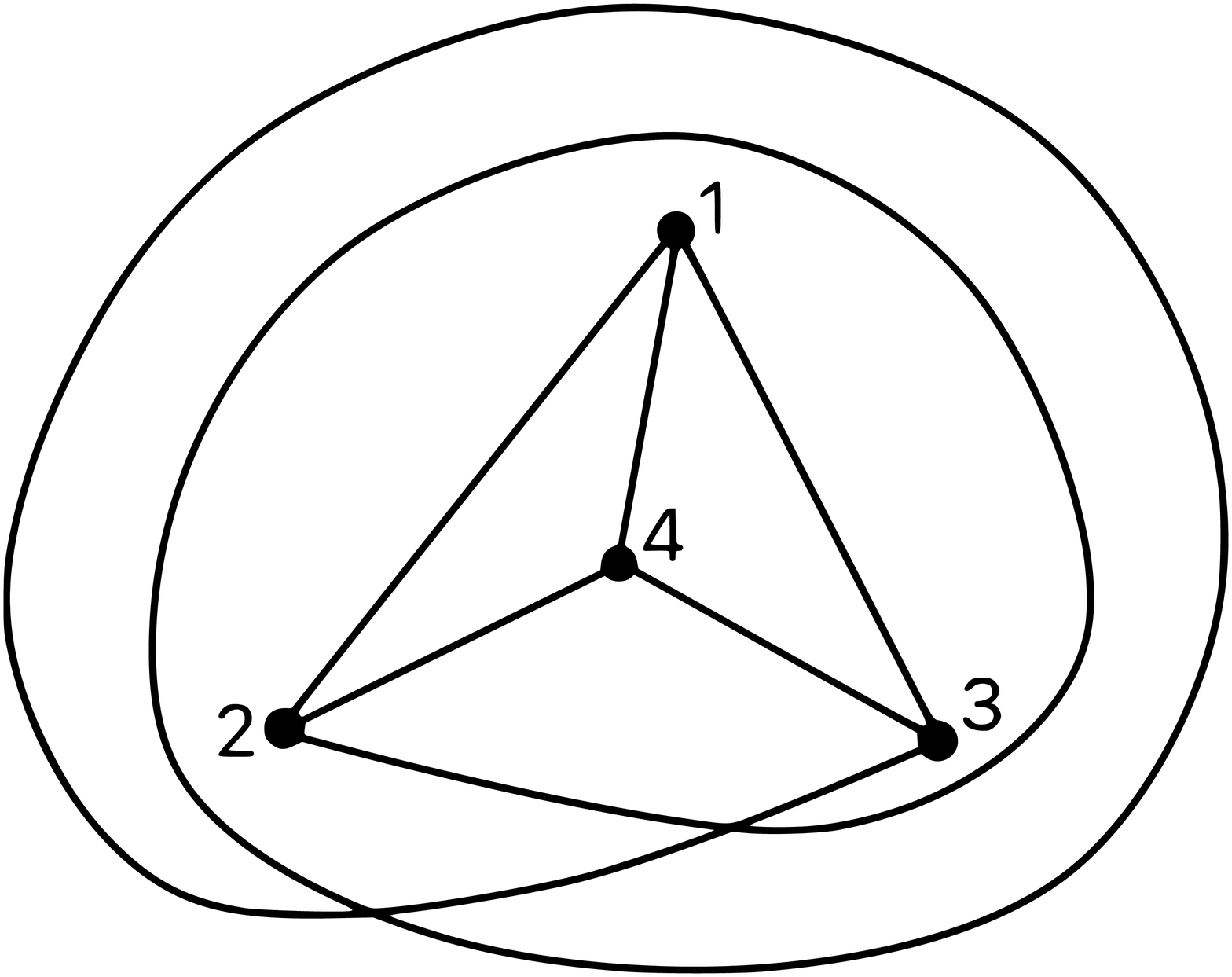}  
\caption{An almost embedding $f:K_4\to\R^2$ such that $w_f(123, 4)=3$}\label{f:123k4}
\end{figure}

\smallskip
\textbf{\ref{p:k3}.} (b) See Figure~\ref{f:123k4} for $n=3$. 
To obtain a map for (a) remove the images of edges issuing out of $f(4)$, and replace $f(4)$ by $O$.

\begin{figure}[ht]\centering
\includegraphics[scale=.8]{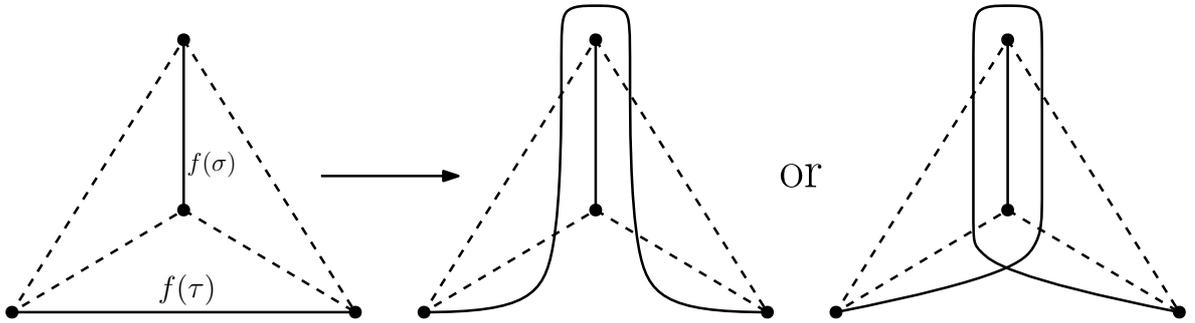}  
\caption{Finger moves (for a map $f$, of an edge $\tau$, w.r.t. the segment $f(\sigma)$) of the first and the second types respectively }\label{f:finger-move}
\end{figure}

\smallskip
\textbf{\ref{e:radonae}.} \emph{Hints.} (a) A regular triangle with its center and edges connecting the center to the vertices form a map $f:K_4 \to \R^2$; assume $f(1),f(2),f(3)$ are the vertices of the triangle numbered counterclockwise.
Then $w_f(C_j, j) = 0$ for every $j \in [3]$, and $w_f(C_4, 4) = 1$.
Make $|n|$ finger moves of the edge $23$ w.r.t. the segment $f(14)$ (Figure~\ref{f:finger-move}, cf. Figure~\ref{f:fingermove}) of the first/second type if $n$ is negative/positive respectively.
The obtained map $f_1$ is as required.

(b) In the next paragraph we construct an almost embedding for $\sum_{j=1}^4 (-1)^j n_j = 1$.
To construct an almost embedding $f$ for $\sum_{j=1}^4 (-1)^j n_j = - 1$ take an almost embedding $g$ for $m_1 = n_2,~ m_2 = n_1,~ m_3 = -n_3,~ m_4 = -n_4$, where $\sum_{j=1}^4 (-1)^j m_j = 1$.
Then $f = g \circ \sigma$, where $\sigma : K_4 \to K_4$  is a permutation that interchanges the vertices $1$ and $2$. 

Consider the map $f_1$ from the proof of (a) for $n = n_1$.
Make $|n_2|$ finger moves of the edge $13$ w.r.t. the segment $f_1(24)$ of the first/second type if $n_2$ is positive/negative respectively.
Denote the obtained map by $f_2$.
Make $|n_3|$ finger moves of the edge $12$ w.r.t. the segment $f_2(34)$ of the first/second type if $n_3$ is negative/positive respectively.
The obtained map $f_3$ is as required.

\begin{figure}[ht]\centering
\includegraphics[scale=1.2]{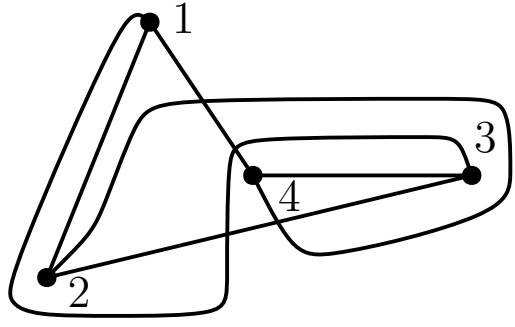}
\caption{An almost embedding $f:K_4\to\R^2$ such that $\sum_{j=1}^4 (-1)^j w_f(C_j, j) = 3$}\label{f:k4-pm3}
\end{figure}

\smallskip
\textbf{\ref{c:radonae}.} \emph{Hint.} (a) See Figure~\ref{f:k4-pm3}.

% (b) Use finger moves of an edge w.r.t. non-self-intersecting polygonal line, defined analogously to the proof of Assertion \ref{e:radonae}.b.  

\smallskip
\textbf{\ref{p:k5-emap}.} (a) No. 
A regular pentagon $f(1)\ldots f(5)$ with all the diagonals, but without the edge $f(4)f(5)$ forms a map $f : K_5-45 \rightarrow \R^2$ such that $w_f(123, 4)-w_f(123, 5) = 0$.

\begin{figure}[ht]\centering
\includegraphics[scale=0.5]{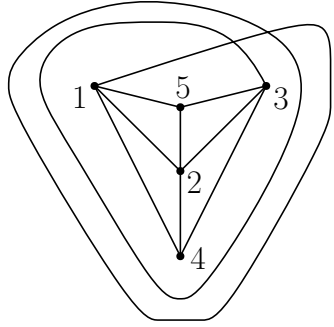}   
\caption{An almost embedding of $K_5$ without an edge $45$ such that $w_f(123, 5) = 3$}\label{f:k5e}
\end{figure}

(b) See Figure~\ref{f:k5e} for $n=3$.

%\newpage
%\section*{Answers, hints and solutions presented at the final discussion}
 
\smallskip
\textbf{\ref{p:triod}.} 
(a) We have 
$$\wu(l_1, l_2, l_3) \stackrel{(1)}{=} 
%2\left(w'(l_2\cup l_3, A_1)+w'(l_1\cup l_2, A_3)+w'(l_3\cup l_1, A_2)\right) \stackrel{(1)}{=}
2\left( \frac{\angle A_2A_1A_3}{2\pi}+k_1 + \frac{\angle A_3A_2A_1}{2\pi}+k_2 + \frac{\angle A_1A_3A_2}{2\pi}+k_3\right)=$$
$$=2(k_1+k_2+k_3)+2\frac{\angle A_2A_1A_3+\angle A_3A_2A_1+\angle A_1A_3A_2}{2\pi}=2(k_1+k_2+k_3)+1$$
for some integers $k_1, k_2, k_3$. 
Here equality (1) follows by (*) and Assertion \ref{p:w'}.a.

\begin{figure}[ht]\centering{
\includegraphics[scale=0.25]{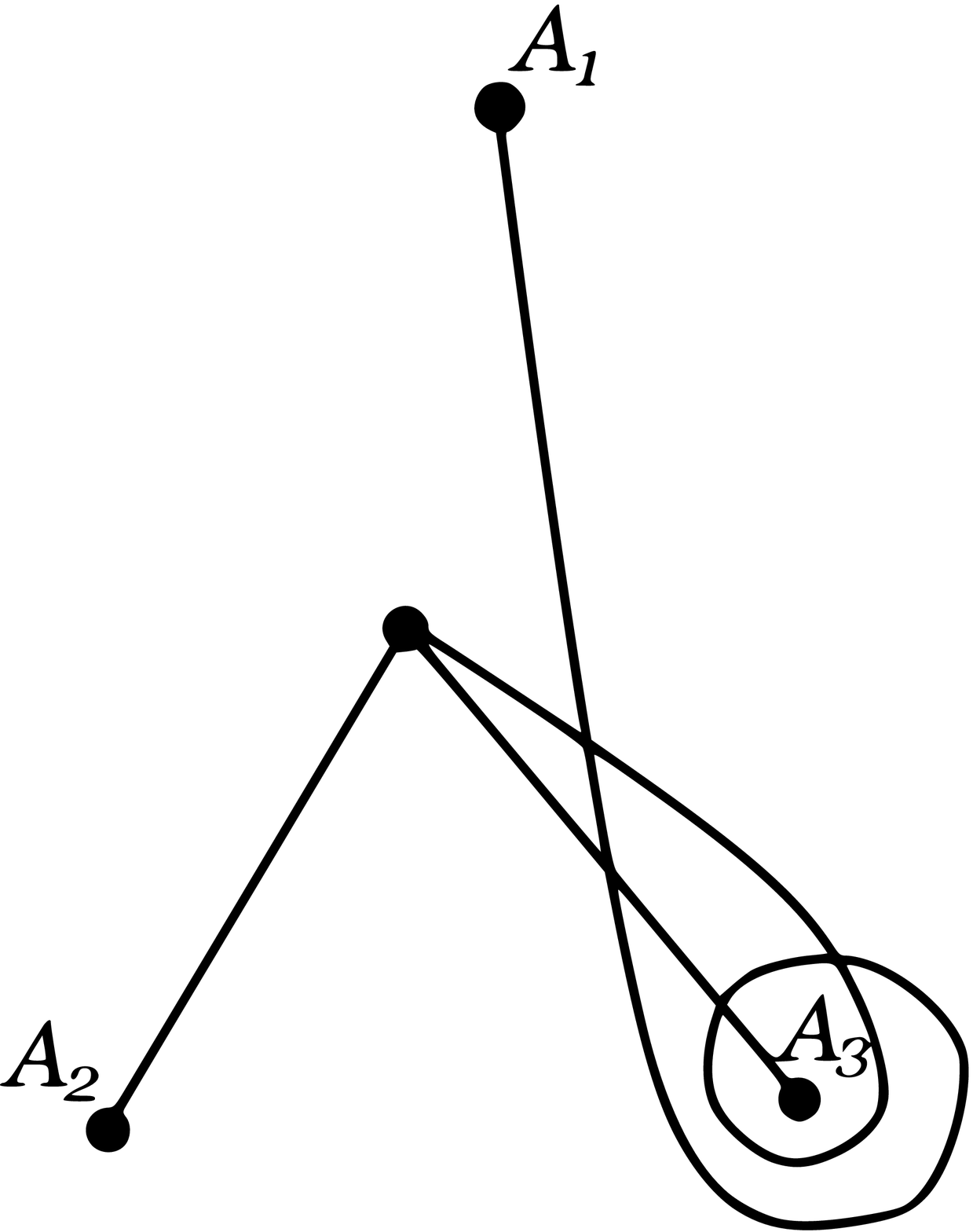} \qquad \includegraphics[scale=0.3]{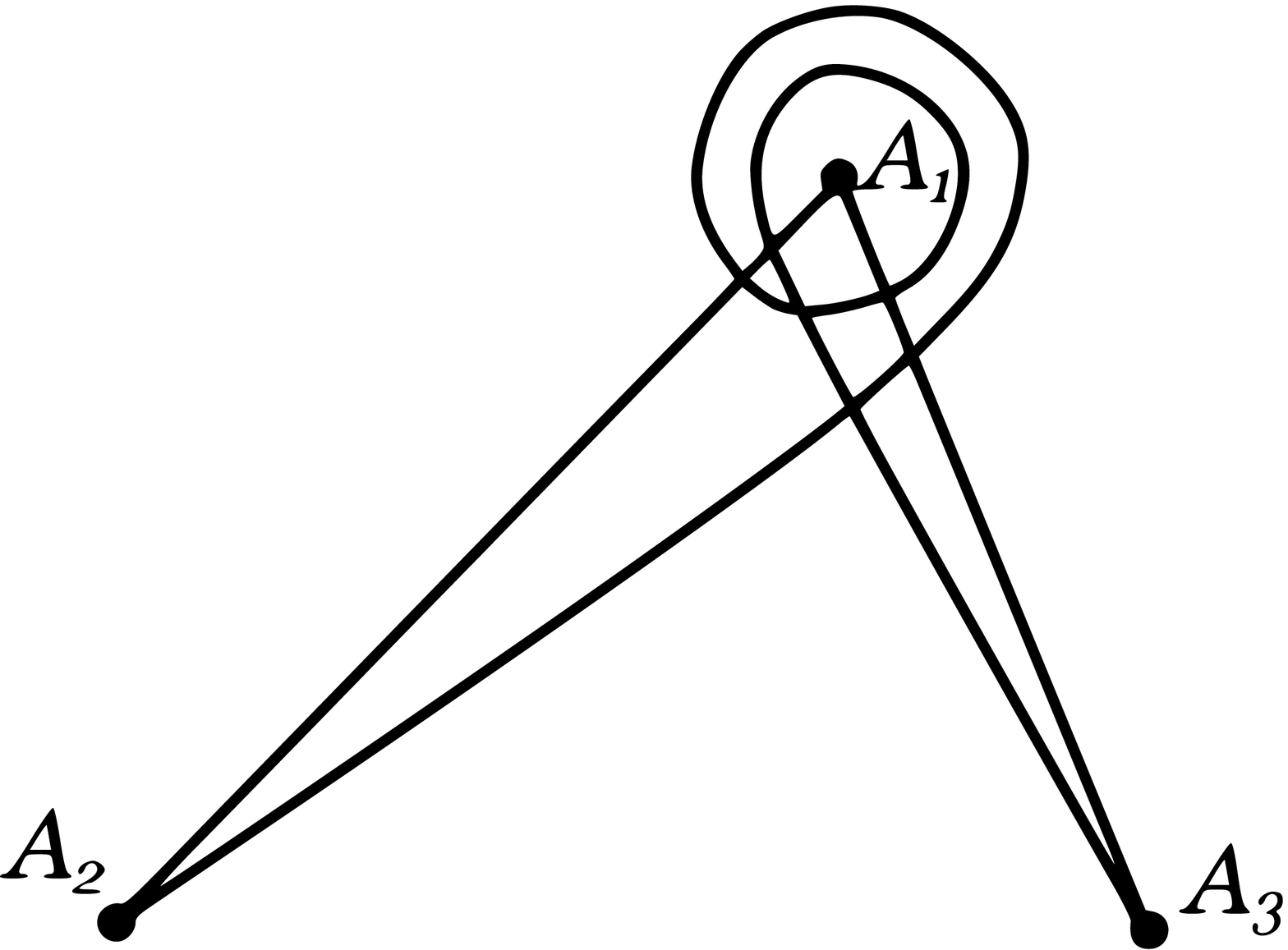}}
\caption{Three polygonal lines whose triodic/cyclic (left/right) Wu number equals $5$}
\label{f:triod5}
\end{figure}

(b) Consider segments $l_1,l_2,l_3$ having a common point $O$ in the interior of triangle $A_1A_2A_3$ such as Figure~\ref{f:triod}. Make $|n|$ finger moves (Figure~\ref{f:fingermove}) of the edge $l_1$ w.r.t. the vertex $A_3$ of the first/second type for negative/positive $n$ respectively. 
Then 
$$\wu(l_1, l_2, l_3)=2\left( \frac{\angle A_2A_1A_3}{2\pi}+\frac{\angle A_3A_2A_1}{2\pi}+\frac{\angle A_1A_3A_2}{2\pi}+n\right)=1+2n.$$ 
For example, in Figure~\ref{f:triod5}, left, 
$$w'(\overline{l_2} l_3, A_1)=\frac{\angle A_2A_1A_3}{2\pi},\quad w'(\overline{l_3} l_1, A_2)=\frac{\angle A_3A_2A_1}{2\pi},\quad w'(\overline{l_1} l_2, A_3)=\frac{\angle A_1A_3A_2}{2\pi}+2,$$ 
so $\wu(l_1, l_2, l_3)=2\cdot 2+1=5$. 

(c) \emph{Hint.} This is proved by induction on the number of segments in $l_1 \cup l_2 \cup l_3$. 

\smallskip
\textbf{\ref{p:off}.} (a) We have 
$$\wu(l_1, l_2, l_3) \stackrel{(1)}{=} 
%2\left(w'(l_2,A_1)+w'(l_1,A_3)+w'(l_3,A_2)\right) \stackrel{(2)}{=}
2\left( \frac{\angle A_2A_1A_3}{2\pi}+k_1+ \frac{\angle A_3A_2A_1}{2\pi}+k_2+ \frac{\angle A_1A_3A_2}{2\pi}+k_3 \right)=$$
$$=2(k_1+k_2+k_3)+2\frac{\angle A_2A_1A_3+\angle A_3A_2A_1+\angle A_1A_3A_2}{2\pi}=2(k_1+k_2+k_3)+1$$ for some integers $k_1, k_2, k_3$. 
Here equality (1) follows by (**) and Assertion \ref{p:w'}.a.
    
(b) Consider triangle $A_1A_2A_3$ with sides $l_2, l_3, l_1$, see Figure~\ref{f:triod}.
Make $|n|$ finger moves (Figure~\ref{f:fingermove}) of the edge $l_2$ w.r.t. the vertex $A_1$ of the first/second type for negative/positive $n$ respectively.
Then $$\wu(l_1, l_2, l_3)=2\left( \frac{\angle A_2A_1A_3}{2\pi}+n+ \frac{\angle A_3A_2A_1}{2\pi}+ \frac{\angle A_1A_3A_2}{2\pi}\right)=2n+1.$$
For example, in Figure~\ref{f:triod5}, right, 
$$w'(l_2, A_1)=\frac{\angle A_2A_1A_3}{2\pi}+2,\quad w'(l_3, A_2)=\frac{\angle A_3A_2A_1}{2\pi}, \quad w'(l_1, A_3)=\frac{\angle A_1A_3A_2}{2\pi},$$ 
so $\wu(l_1, l_2, l_3)=2\cdot 2+1=5$.

%AS where is this used? In order to obtain an almost embedding from an almost embedding by a finger move (Figure \ref{f:fingermove}), it is necessary that the edge $\tau$ is not adjacent to any edge issuing out of the vertex $a$.

\begin{figure}[ht]
\centerline{\includegraphics[width=4.5cm]{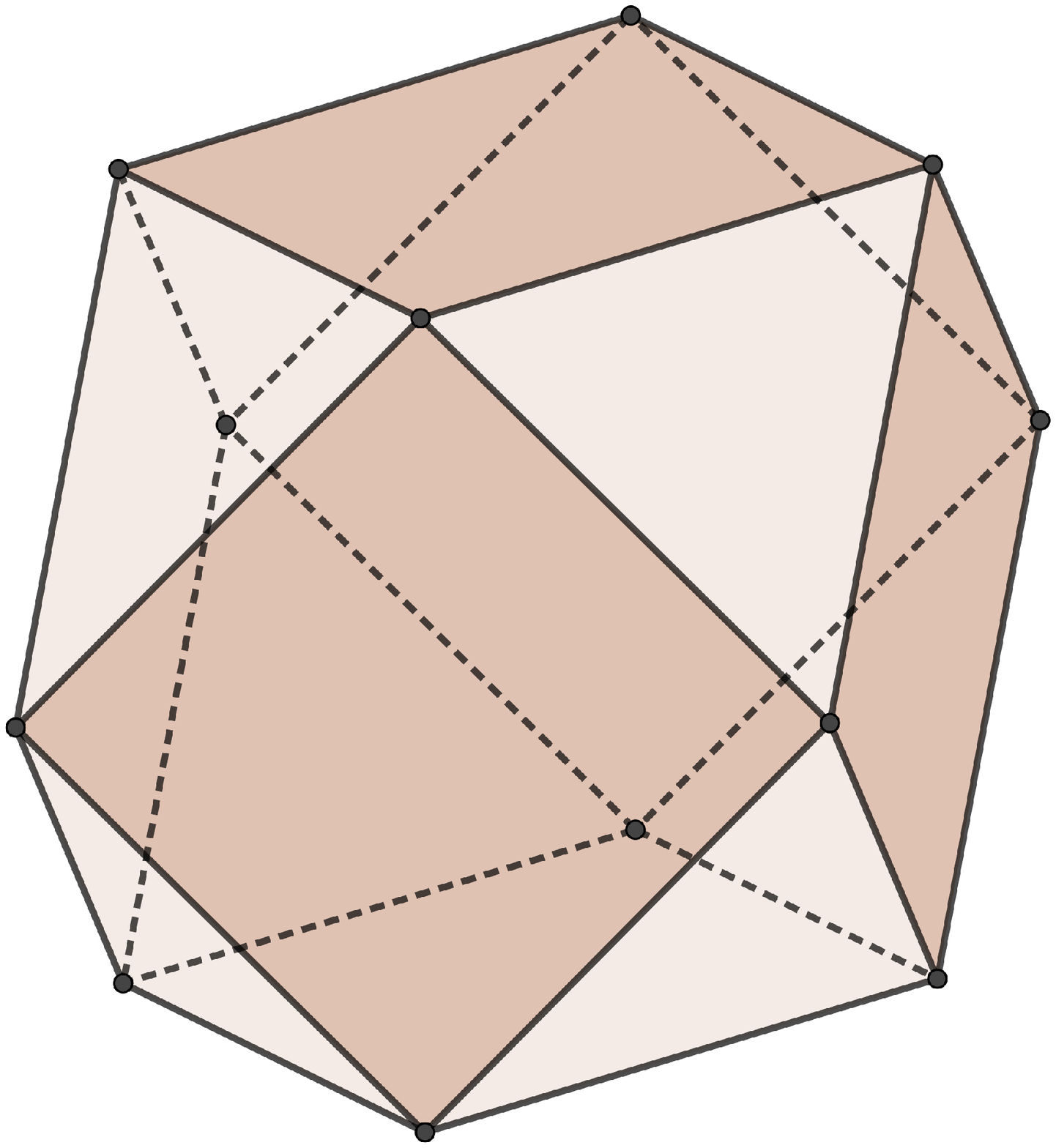}
\qquad\includegraphics[width=6cm]{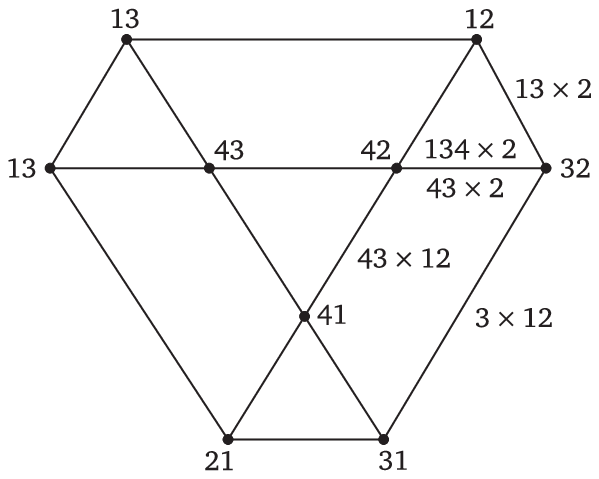}}
\caption{Left: a magic cuboctahedron. 
%;  the complement to triangular faces is $K_4^{\underline2}$. \newline
Right: the magic:
%same with some explanations; 
the figure does not show the invisible part whose projection is obtained from the pictured projection by rotation through $\pi/3$; the lower 13 should be 23; $4\times 123$ is the central triangle, $123\times4$ is the invisible central triangle,  $24\times 13$ is the bottom left trapezoid, $13\times24$ is the upper right invisible trapezoid,
$\text{off}(123)=12, 13, 23, 21, 31, 32$ is the outer cycle of length 6, $\text{triod}(123, 4)= 12, 14, 13, 43, 23, 24, \ldots,$ where dots denote the part symmetric to the written part (obtained replacing $xy$ by $yx$) 
%is certain cycle of length 12. 
(each of the cycles $\text{off}(123)$ and $\text{triod}(123, 4)$ splits the cuboctahedron into two equal parts)}
\label{f:del3}
\end{figure}

\smallskip
\textbf{\ref{p:wu-conj}.} 
{\it Hint:} see Figure \ref{f:del3} and \cite[\S2, \S5]{ADN+}; $\text{off}(123) + \text{triod}(123,4) =$ 
\linebreak
$123\times 4 + 4\times 123 + \partial(13\times 24) + \partial(24\times 13) + \partial(12\times 34) + \partial(34\times 12) + \partial(14\times 23) + \partial(23\times 14)$. 

%This appendix is intended for participants successfully solving other problems.   

%\newpage

{\it Books, surveys and expository papers in this list are marked by the stars.}

\end{document}